\documentclass[11pt]{article}
\usepackage{graphicx}
\usepackage{amssymb}
\usepackage{epstopdf}
\usepackage{amstext,amsmath}

\usepackage{array}
\input xy
\xyoption{all}

\def\GL{\mbox{GL}}

\def\R{\mathbb{R}}
\def\C{\mathbb{C}}
\def\Q{\mathbb{Q}}
 \def\Z{\mathbb Z}
 
\def\Mod{\mbox{Mod}}

\def\sF{\mathcal{F}}
\def\sC{\mathcal C}
\def\sW{\mathcal W}
\def\sR{\mathcal R}
\def \sS{\mathcal S}
\def \sA{\mathcal A}

\def\sV{\mathcal V}
\def\sE{\mathcal E}

\def\sO{\mbox{O}}

\def\s{\mathfrak{s}}
\def\oo{\mathfrak{o}}

 \def\qed{\hfill\framebox(5,5){}}
\DeclareMathAlphabet{\mathpzc}{OT1}{pzc}{m}{it}
\def\sign{\mathfrak{s}}
\def\dil{d}

\def\f{\mathfrak f}
\renewcommand{\bold}[1]{\smallskip \noindent {\bf \boldmath #1 }\nopagebreak[4]}

\newtheorem{theorem}{Theorem}[section]
\newtheorem{corollary}[theorem]{Corollary}
\newtheorem{lemma}[theorem]{Lemma}

\newtheorem{problem}[theorem]{Problem}
\newtheorem{remark}[theorem]{Remark}
\newtheorem{example}[theorem]{Example}
\newtheorem{question}[theorem]{Question}
\newtheorem{proposition}[theorem]{Proposition}
\newtheorem{observation}[theorem]{Observation}

\title{Mapping classes associated to mixed-sign Coxeter graphs}
\author{Eriko Hironaka\footnote{This work was partially supported by a grant from the Simons Foundation (\#209171 to Eriko Hironaka).}} 
\begin{document}
\maketitle

\begin{abstract}
We define a generalization of Coxeter graphs and an associated Coxeter system and Coxeter mapping class.  These
can be used to construct periodic Coxeter mapping classes on surfaces with arbitrarily large genus, preserving lots
of symmetries.   The periodic mapping classes can in turn be used to construct sequences of pseudo-Anosov 
mapping classes with bounded normalized dilatation and arbitrarily high genus.   We show that the smallest known 
accumulation point of normalized dilatations can be realized by such a sequence.
\end{abstract}


\section{Introduction}

Let $S$ be a compact oriented surface of finite type.  A mapping class  on $S$ is an isotopy class
of self-homeomorphisms of $S$.  Of particular interest are the {\it irreducible} 
mapping classes, which do not fix any nontrivial
multi-curves on $S$.  Such mapping classes $\phi$ are  called 
{\it pseudo-Anosov} 
and have the property that for some pair of $\phi$-invariant, singular foliations $\mathcal F^+$ and $\mathcal F^-$,
with transverse measures $\nu^+$ and $\nu^-$, 
we have $\nu^\pm(\phi(x),\phi(y)) = \lambda^{\pm 1} (x,y)$, for $x$ and $y$
points on the nonsingular locus of $\sF^\pm$ (\cite{Thurston88}, \cite{FLP79}).   When the leaves
in $\mathcal F^\pm$ can be given a consistent orientation, 
we say that $\phi$ is an {\it orientable pseudo-Anosov mapping class}.    This is equivalent to the statement that
the dilatation $\lambda(\phi)$ equals the homological dilatation of $\phi$, that is, the spectral radius of the action
of $\phi$ on the first homology of $S$.
The minimum possible value of $\lambda(\phi)$ for a fixed surface $S$ is known
only for surfaces with small genus and number of boundary components \cite{SoKoLos02} \cite{HamSong05} \cite{CH08} \cite{LT11}.  
Slightly more is known for the case of orientable pseudo-Anosov mapping classes \cite{HK:braidbounds} \cite{LT09}.

The minimum dilatation of pseudo-Anosov
mapping classes on a surface $S$ is known to go to 1 as the absolute value of the topological Euler characteristic
$|\chi(S)|$ goes to infinity.  The first proof of this was given in \cite{Penner91}, and the orientable case was shown
in \cite{HK:braidbounds}.  The minimum value of $L(S,\phi)$ for $S = S_{g,n}$ is bounded for $(g,n)$ on any line of rational slope
through the origin  \cite{Penner91} \cite{HK:braidbounds} \cite{Tsai08} \cite{Valdivia:gn}, but unbounded on any line with
$g$ constant: $g = g_0 \ge 2$ \cite{Tsai08}.      We will say that a
pseudo-Anosov mapping class has {\it small dilatation} if the 
 {\it normalized dilatation} 
 $$
L(S,\phi)= \lambda(\phi)^{|\chi(S)|}
$$
is small.   

\begin{question}\label{smalldil-ques} What do small dilatation pseudo-Anosov mapping classes look like?
\end{question}

In this paper, we approach question~\ref{smalldil-ques} by studying a simple construction of pseudo-Anosov
mapping classes from generalized Coxeter graphs.    We call such mapping classes {\it mixed-sign Coxeter
mapping classes}.  Our main result is the following.

\begin{theorem}\label{main1-thm} The set of normalized dilatations of 
 pseudo-Anosov mixed-sign Coxeter mapping classes has accumulation point $\gamma_0^4$.
 \end{theorem}

The minimum accumulation point  $\ell_0$  of $L(S,\phi)$ over all pseudo-Anosov mapping classes satisfies
\begin{eqnarray}\label{bound-eqn}
2 < \ell_0 \leq  \left (\frac{3 + \sqrt{5}}{2} \right )^2 = \gamma_0^4 \approx 6.8541,
\end{eqnarray}
where $\gamma_0$ is the golden mean (\cite{Hironaka:LT} \cite{AD10} \cite{KT11}).   The lower bound  comes from properties
of directed graphs (see \cite{Penner91}
\cite{McMullen:Poly}).
The upper bound $\gamma_0^4$ in (\ref{bound-eqn}) is achieved by the {\it simplest hyperbolic braid}
(see, for example, \cite{McMullen:Poly} \cite{Hironaka:LT}).  
Smaller normalized dilatations exist, for example, for the once-punctured torus, the smallest normalized
dilatation is $\frac{3 + \sqrt{5}}{2}  \approx 2.61803$, and the smallest normalized dilatation for
a closed genus $2$ surface $|x^4 - 2x^2 - 2x +1| \approx 5.27451$.

\subsection{Fibered face theory}
By  Thurston's fibered face theory \cite{Thurston:norm}
and a result of Fried \cite{Fried82}, the rational points on fibered faces of hyperbolic 3-manifolds
parameterize all pseudo-Anosov mapping classes on compact oriented surfaces (this is explained
in more detail below.
Since the simplest hyperbolic braid has normalized dilatation equal to $\gamma_0^4$, and its
fibered face has positive dimension,
 the values of $L$ are dense in the interval $[\gamma_0^4,\infty)$
(see  \cite{Hironaka:LT}).   The genus 1 and genus 2 examples mentioned above
lie on fibered faces consisting of a single point, and can be considered as isolated points in the
space of pseudo-Anosov mapping classes.

    The Universal Finiteness Theorem of Farb, Leininger and Marglit \cite{FLM09}
implies that there is a finite collection of mapping tori (up to homeomorphism) for
pseudo-Anosov mapping classes with no interior singularities whose normalized dilatations is bounded. 
Fibered face theory implies that the dynamical structure of small dilatation mapping classes are governed by a finite set 
of ``templates" defined by pseudo-Anosov flows on hyperbolic 3-manifolds \cite{Thurston:norm} \cite{Fried82} \cite{FLM09}.

More precisely,  if a hyperbolic $3$-manifold $M$ is fibered and has first Betti number $b_1(M)$ greater than or equal to 2, 
then there is at least one open cone, called a {\it fibered cone}, in $\R^{b_1(M)}$ whose primitive integral points parameterize
 fibrations of $M$ to $S^1$.     This fibered cone is of the form $F \cdot \R^+$, where $F$ is a fibered face
 of the Thurston norm ball for $M$ in $H^1(M;\R)$ (see \cite{Thurston:norm} and Section~\ref{fiberedface-sec} for
 definitions).   By a result of Fried \cite{Fried82}
the normalized dilatation function $L$ extends to a continuous convex function on the cone $F \cdot \R^+$ so that 
$\log L$ is homogeneous
of degree $-1$.   It follows that on every compact subset of $F$ there are
pseudo-Anosov mapping classes with unbounded topological Euler characteristic and bounded normalized dilatation.
Conversely, if we restrict to pseudo-Anosov mapping classes with no interior singularities, 
all such infinite families lie on a finite union of such fibered faces by the Universal Finiteness Theorem.   Thus,
Question~\ref{smalldil-ques} is strongly connected to the question.

\begin{question} What do mapping classes associated to rational points on a single fibered face look like?
\end{question}

So far there has been no explicit characterization of which 3-manifolds realize minimum normalized diatations
smaller than a given bound, nor has there been a simple constructive description of the infinite family of pseudo-Anosov 
maps associated to the rational points on an arbitrary fibered face.

Our main new technique in this paper is to show that there are natural sequences of mapping classes 
associated to generalized Coxeter graphs that correspond to Cauchy sequences on fibered 
faces.  We show that generalized Coxeter mapping classes contain a large family of
periodic mapping classes on surfaces with high symmetry.  These provide useful building blocks for defining small 
dilatation pseudo-Anosov mapping classes.   We observe that the
 minimum dilatation orientable pseudo-Anosov mapping classes for genus $g = 2,3,4,5$ found by Lanneau and
 Thiffeault \cite{LT09} can be defined from generalized Coxeter graphs (Table~\ref{smallorientable-table}).  We also 
 define a sequence $(S_n,\phi_n)$
of pseudo-Anosov mapping classes whose normalized dilatations converge to $\gamma_0^4$ (Section~\ref{twistexample-sec}).  

\subsection{Mixed-sign Coxeter reflection groups.}  
A mixed-sign Coxeter system is a generalization of  classical Coxeter systems.  A
{\it (simply-laced, classical) Coxeter graph} is a finite connected
graph with no self- or double-edges.  A Coxeter graph $\Gamma$ defines a
{\it Coxeter reflection group}
$\sW_\Gamma \subset \GL(\R^\sV)$, where $\sV$ is the set of vertices of $\Gamma$, and $\R^{\sV}$ is the vector space of
$\R$ valued functions on $\sV$.  
The group $\sW_\Gamma$ is generated by a distinguished set of reflections
$S = \{s_1,\dots,s_k\}$ preserving the bilinear form $B_\Gamma = 2I - A_\Gamma$,
where $A_\Gamma$ is the adjacency matrix of $\Gamma$, and $I$ is the identity matrix (see  \cite{Humphreys:Coxeter}
and Section~ ).
The {\it Coxeter element} $\omega_\Gamma = s_1 \cdots s_k$ can be used to classify the Coxeter system.
For example, if $\Gamma$ is connected, then $B_\gamma$ defines a {\it spherical}, {\it affine}, or {\it higher rank} 
geometry on $\R^\sV$ if and only if
$\omega_\Gamma$ is finite order, has infinite order but spectral radius 1, or has spectral radius greater than 1,
respectively (see \cite{ACampo:Coxeter}).   The pair $(\R^\sV,B_\Gamma)$ and 
reflection group $\sW_\Gamma$ with its distinguished generating reflections
and Coxeter element is called the {\it Coxeter system} associated to $\Gamma$.

%

A {\it mixed-sign Coxeter graph} is a Coxeter graph with labels 
$$
\s : \sV \rightarrow \{1,-1\}.
$$
The associated reflection group $\sW_{\Gamma,\s} \subset \GL(\sR^\sV)$ is defined in an analogous way as the classical
Coxeter reflection group, except that the bilinear form is replaced by
$B_{\Gamma,\s} = 2I_\s - A_\Gamma$, where $I_\s$ is the diagonal matrix with entries $\s(v_i)$ on the diagonal.  
The relation between properties of $\omega_{\Gamma,\s}$ and properties of
$\sW_{\Gamma,\s}$  is more subtle  than in the classical case (see  \cite{Armstrong:thesis}),
and, for example, the finite order mixed-sign Coxeter reflection groups are not yet classified.
However,  a useful formula for the Coxeter element in terms of the adjacency matrix still holds
in this setting (see Proposition \ref{Howlett-prop}).

\subsection{Constructing small dilatation pseudo-Anosov mapping classes.}  The first construction of infinite sequences of pseudo-Anosov
mapping classes $(S_n,\phi_n)$, where $|\chi(S_n)|$ is unbounded and $L(S_n,\phi_n)$ is bounded, was given by Penner
in \cite{Penner91} (see also \cite{Bauer:thesis} \cite{Valdivia:gn}).    There is no evidence, however, that
the conjectural minimum accumulation point for normalized dilatations $\gamma_0^4$ can be achieved by a Penner-type sequence.

 Given a mapping class $(S,\phi)$, 
    the {\it closure} of $(S,\phi)$ is the pair $(\overline S,\overline \phi)$ where $\overline S$ is the closed surface obtained from $S$
 by filling in boundary components with disks, and $\overline \phi$ is the isotopy class of  the extension to $\overline S$ of any homeomorphism
 in the equivalence class of $\phi$.

  \begin{question} For which closed oriented surfaces $S$ can the minimum dilatation pseudo-Anosov mapping classes 
  on $S$  be realized as the closure of a 
  mixed-sign Coxeter mapping class?  \end{question}
  
  The minimum dilatations are unknown for genus $g \ge 3$, however, in the orientable case, when 
  $\lambda_{\mbox{hom}}(\phi) = \lambda(\phi)$, the smallest dilatations, together with explicit definitions of their
  monodromy, were found for genus $g = 2,3,4,5$
   by E. Lanneau and  J.-L. Thiffeault \cite{LT09}.  We make the following observation (see Table~\ref{smallorientable-table}).
  
\begin{observation}\label{main1-obs} The minimum dilatations for orientable mapping classes for genus 2,3,4, and 5 are realizable as
the closures of  mixed-sign Coxeter mapping classes.
\end{observation}

\begin{remark} {\em The minimum dilatation for orientable mapping classes for genus 7 was found in  \cite{AD10} and in \cite{KT11}, and 
for genus 8 in \cite{Hironaka:LT}.  These examples were found using fibered face theory, and there is as yet no
intrinsic description.}
\end{remark}


\subsection{Lanneau-Thiffeault Question.}
For orientable pseudo-Anosov mapping classes Lanneau and Thiffeault ask (see \cite{LT09}) whether the minimum dilatation
for orientable pseudo-Anosov mapping classes is the largest root $\lambda_g$ of polynomials of the form
$$
LT(x) = x^{2g} - x^{g+1}-x^g -x^{g-1} + 1,
$$
for all even $g$.  In \cite{Hironaka:LT} it is shown that  for even $g$, there are orientable genus $g$ pseudo-Anosov mapping classes
$(S_g,\phi_g)$ with exactly two singularities whose dilatation equals $\lambda_g$.  This sequence of mapping classes corresponds to a convergent
sequence of points on a single fibered face, and we have
$$
\lim_{g \rightarrow \infty} L(S_g,\phi_g) = \lim_{g \rightarrow \infty} (\lambda_g)^{2g} = \gamma_0^4.
$$
Thus, the smallest known accumulation point for normalized dilatations, is also achieved by a sequence of orientable
pseudo-Anosov mapping classes implying the following.

\begin{theorem}\label{main2-thm} The set of normalized dilatations of 
closures of orientable pseudo-Anosov mixed-sign Coxeter mapping classes has accumulation point $\gamma_0^4$.
\end{theorem}

Theorem~\ref{main1-thm} follows from Theorem~\ref{main2-thm}.

To find graphs associated to small dilatation mapping classes, we use the heuristic principle that a mapping class
with small dilatation should be a slight perturbation of a mapping class that is periodic.    Roughly speaking,
the Coxeter element of a Coxeter graph is conjugate to the homological action of its corresponding mapping
class.  For classical Coxeter systems, the Coxeter element is periodic if and only if the Coxeter system is
spherical.  These are classified, and furthermore, the spectral radius of Coxeter elements is monotone with
respect to graph inclusion.  Thus, for mapping classes associated to classical Coxeter systems, there is
a universal lower bound on the dilatations of associated pseudo-Anosov maps (see \cite{McMullen:Coxeter} and
\cite{Leininger04}).   

For mixed-sign Coxeter graphs, the spectral radius of Coxeter elements is not monotone.  Thus, it is possible
to find large and complicated graphs whose Coxeter elements are periodic.  These are useful tools for building
pseudo-Anosov mapping classes with large complexity (as measured by the Euler characteristic of the surface) 
and small dilatation.   More precisely, we find sequences of spherical mixed-sign graphs $K_n$ with 
increasing number of vertices, a fixed graph $\Gamma_0$ and concatenations $\Gamma_n = K \sharp K_n$, so that
the spectral radius $\lambda(\Gamma_n)$ of the Coxeter element of $\Gamma_n$ goes to 1 quickly, i.e.,
$$
\log \lambda(\Gamma_n) \asymp \frac{1}{n}.
$$

Murasugi sums provide a way to translate between graph theoretic perturbations, and perturbations of mapping classes.
Using this idea, we prove that the sequences of mapping classes associated to $K_n$ are the Murasugi sum
of simple to understand mapping classes.  We further show  that mapping classes $(\Sigma_n,\rho_n)$ associated to $\Gamma_n$
correspond are fibrations of a single 3-manifold.  Thus the Murasugi sum of $(\Sigma_n,\rho_n)$ with any
fixed mapping class defines a sequence $(S_n,\phi_n)$ whose mapping tori are homeomorphic.  
We call such mapping classes $(S_n,\phi_n)$
twisted mapping classes. Fibered face theory then implies that
the sequence of normalized dilatations has accumulation points either in the interior or on the boundary of a fibered face.
By arranging for the former, we obtain many sequences of mapping classes with bounded normalized dilatations.

In particular, we show that for our particular choice $\Gamma_n$, the mapping classes
$(S_{\Gamma_n},\phi_{\Gamma_n})$ correspond  to points on the fibered face of the
 $6{}_2^2$-link complement and converge to the minimum $(S_0,\phi_0)$ for normalized dilatation.  
 By continuity of normalized dilatation on fibered faces, it then follows that 
 $$
 \lim_{n \rightarrow \infty} L(S_{\Gamma_n},\phi_{\Gamma_n}) = \gamma_0^4.
 $$
 
 \subsection{Organization.}  We begin in Section~\ref{alg-sec} with definitions of mixed-sign Coxeter systems $\Gamma,\s$ and their
 properties, particularly in relation to the Artin group $\sA_\Gamma$ associated to $\Gamma$.
We also define a surface $S_\Gamma$ associated to $\Gamma$ given an ordering on the vertices
 and a fat graph structure.   Then the Artin group $\sA_\Gamma$ acts as mapping classes on $S_\Gamma$, and 
 for each choice of signs $\s$ on $\Gamma$, we define a particular element $\phi_{\Gamma,\s}$ depending on an ordering
 on $\Gamma$.  If $\Gamma$ is bipartite, then this is the same as examples studied by Thurston \cite{Thurston88} and
 Leininger \cite{Leininger04}.   In Section~\ref{twist-sec}, we define twisted mapping classes and prove
 Theorem~\ref{main2-thm}.

\bold{Acknowledgments:}  I am grateful to the Tokyo Institute of Technology and University of
Tokyo for their support during the writing of this paper, and J. F. Valdez and J. Mangahas for helpful conversations.

\section{Mixed-sign Coxeter systems and associated mapping classes.}\label{alg-sec}

In this section we define a correspondence between various objects associated to a mixed-sign Coxeter graph.
The following table summarizes the main objects.
\begin{center}
\begin{tabular}{|l|l|l|}
\hline
{\bf Data} & {\bf Object} & {\bf Notation}\\
\hline
mixed sign  graph $\Gamma,\s$ &Mixed-sign Coxeter system &$(\sW_{\Gamma,\s},\sR_{\Gamma,\s})$\\
\hline
 ordered mixed sign graph $\Gamma,\s$ &Mixed-sign Coxeter element &$\omega_{\Gamma,\s}$\\ 
 \hline
graph $\Gamma$ &Artin group and representation& $\rho_{\Gamma} : \sA_{\Gamma} \rightarrow GL(\R^{\sV})$\\
 \hline
ordered  fatgraph  $\Gamma$ & geometric realization &$(S_\Gamma,\sC_\Gamma)$\\
\hline
ordered mixed sign fatgraph $\Gamma,\s$ &Coxeter mapping class & $(S_\Gamma,\phi_{\Gamma,\s})$\\
\hline
\end{tabular}\label{objects-table}
\end{center}

Given an ordered mixed-sign Coxeter  graph,  $(\Gamma,\s)$  we define a reflection group
$\sW_{\Gamma,\s}$ in Section~\ref{Coxeter-sec},  and a representation of the Artin group $\sA_\Gamma$ 
$$
\rho_{\Gamma} : \sA_{\Gamma} \rightarrow GL(\R^{\sV})
$$
in Section~\ref{Artin-sec}.  When $\Gamma$
has a fatgraph structure, we define an associated geometric realization $(S_\Gamma,\sC_\Gamma)$ of $\Gamma$,
and a mapping class $\phi_{\Gamma,\s} : S_{\Gamma} \rightarrow S_\Gamma$ (Section~\ref{geo-sec}). 
The groups $\sA_{\Gamma}$ and $\sW_{\Gamma,\s}$ come with $k$ standard generators $\{\sigma_1,\dots,\sigma_k\}$
and $\{s_1,\dots,s_k\}$, respectively.
Let $\sigma$ and $w$ be the 
epimorphisms of the free group $\langle x_1,\dots,x_k\rangle$ to $\sA_{\Gamma}$ and $\sW_{\Gamma,\s}$, where $\sigma(x_i) = \sigma_i$
and $w(x_i) = s_i$.
$$
\xymatrix{
&\langle x_1,\dots,x_k\rangle \ar[r]^{\quad \sigma}\ar[dr]^{w}&\sA_{\Gamma} \ar[r]^{\rho_{\Gamma}\quad}& GL(\R^\sV) \ar[r]^{\eta\quad} &GL(H_1(S_\Gamma;\R))\\
&&\sW_{\Gamma,\s} \ar[r]^{\subset\quad} & GL(\R^{\sV}).
}
$$ 
Let $\omega_{\Gamma,\s} = w (x_1 \cdots x_k)$, called the {\it (mixed-sign) Coxeter element}.  Write $\s(i) = \s(v_i)$, and let 
$\sigma_{\Gamma,\s} = \rho_{\Gamma}(\sigma_1^{\s(1)} \cdots \sigma_k^{\s(k)})$.
The representations $\rho_{\Gamma}$ and $\eta$  preserve symplectic forms, while the elements of $\sW_{\Gamma,\s}$ preserve a
symmetric one, 
but we will show
$$
\omega_{\Gamma,\s} = - \sigma_{\Gamma,\s}
$$
as elements of $GL(\R^{\sV})$ (Proposition~\ref{ArtinCoxeter-prop}).
The element
$\sigma_{\Gamma,\s}$ satisfies
$$
\eta (\rho_{\Gamma}(\sigma_{\Gamma,\s})) = (\phi_{\Gamma,\s})_*
$$
and hence the homological dilatation of $\phi_{\Gamma,\s}$ and the spectral radius of the Coxeter element $\omega_{\Gamma,\s}$
satisfy
$$
\lambda_{\mbox{hom}}(\phi_{\Gamma,\s}) = |\omega_{\Gamma,\s}|
$$
(see Proposition~\ref{hom-prop}).
From this we derive a sufficient condition for a mixed-sign Coxeter mapping class to be pseudo-Anosov  (see 
Proposition~\ref{pAcrit-prop}).

Mixed-sign Coxeter mapping classes are defined on surfaces with boundary, but in some cases the dynamical information contained
in the mapping class extends to the closure of the surface obtained by filling in disks.  This is discussed in Section~\ref{closure-sec}.
The special case of mapping classes associated to bipartite Coxeter graphs is treated in Section~\ref{bipartite-sec}. 
Section~\ref{orientable-sec} gives a table of mixed-sign Coxeter graphs associated to minimum dilatation
orientable pseudo-Anosov mapping classes for small genus.

\begin{remark} \em{ The study of mapping classes using Coxeter graphs and associated reflection
 groups   has a long history
 in algebraic geometry dating back to the 19th century, particularly in the study of complex
 surface singularities
  (see \cite{Saito98}, \cite{Dolgachev08} and references therein). 
The focus in geometric topology has on the other hand been on  representations
 of the Artin group of a Coxeter graph into the mapping class group (see, for example, \cite{LabruereParis01}).  The difference
 comes from the fact the associated 
  bilinear forms  left invariant by the respective automorphism groups are different:  one
 being symmetric and the other symplectic.  From this point
 of view our results concerning special elements of the
  Artin group, and Coxeter elements of the Coxeter reflection 
	 group are part of an overlap in the two theories (see also \cite{Arnold:Singularities} \cite{Hirz:Sing}).}
 \end{remark}

\subsection{Mixed-sign Coxeter reflection group}\label{Coxeter-sec}

In this section we define mixed-sign Coxeter systems as a slight generalization of classical simply-laced Coxeter systems.
Let $(\Gamma,\s)$ be a mixed-sign Coxeter graph
with ordered vertices $\sV=\{v_1,\dots,v_n\}$ and map $\s : \sV \rightarrow \{\pm 1\}$.
For $i,j=1,\dots,k$, define
$$
m_{i,j} = 
\left 
\{
\begin{array}{rl}
1 & \qquad\mbox{if $i=j$}\\
2 & \qquad \mbox{if $v_i$ and $v_j$ are not connected by an edge}\\
3 & \qquad \mbox{if $v_i$ and $v_j$ are connected by an edge.}
\end{array}
\right .
$$

\begin{remark}{\em In the usual definition of Coxeter and Artin groups, the  $m_{i,j}$ are allowed
to vary (for $i \neq j$)  in the set $\{2,3,\dots,\infty\}$.   In this paper, we restrict only to the {\it simply-laced
case} described above.}
\end{remark}

Let $I_\sign$ be the $n \times n$ matrix with $0$ on the off-diagonal, and diagonal entries equal to 
$\s(1),\dots,\s(n)$.  Let $A = [a_{i,j}]$ be the adjacency matrix of $\Gamma$, $A^+$ the upper
triangular part of $A$, and
 $U= I_\sign - A^+$.     Then $U$ 
determines a symmetric form
$$
B = U + U^T
$$
on $V$.

The Coxeter system $(\sW,\sR) = (\sW_{\Gamma,\s},\sR_{\Gamma,\s})$ associated to $(\Gamma,\s)$ is the subgroup of $\GL(\R^{\sV})$
generated by $\sR_{\Gamma,\s} = \{s_1,\dots,s_n\}$, where the $s_i$ are defined by
\begin{eqnarray*}
s_i(v_j) &=& v_j - 2 \frac{B(v_i,v_j)}{B(v_I,v_I)}v_i\\
&=& 
\left \{
\begin{array}{rl}
-v_j &\qquad\mbox{if $i=j$,}\\
v_j + \s(i) a_{i,j} v_i &\qquad\mbox{if $i \neq j$.}
\end{array}
\right .
\end{eqnarray*}
Each $s_i$ can be interpreted as a reflection through the hyperplane perpendicular
to $v_i$ in $\R^{\sV}$ with respect to the symmetric form $B$.
If $\sign \equiv 1$, then $(\sW,\sR)$ is the classical
Coxeter system.

Given a mixed-sign ordered Coxeter graph $(\Gamma,\s)$, we define the {\it Coxeter element} to be
$$
\omega_{\Gamma,\sign} = s_1 \cdots s_n \in \sW_{\Gamma,\sign}.
$$

Coxeter systems $(\sW,\sR)$, and more generally reflection systems, are classified by the type of
associated bilinear form $B$.  If $B$ is positive or negative definite, we say $(\sW,\sR)$ is 
{\it spherical}, if $B$ is positive or negative semi-definite, then we say $(\sW,\sR)$ is
{\it affine}, and if the signature of $B$ is $(p,q)$, where $p$ and $q$ are nonzero we say $(\sW,\sR)$ is {\it higher rank}.  
(see \cite{Humphreys:Coxeter}). For classical Coxeter systems, these are the only cases that occur.

The classical (simply-laced) Coxeter systems $(\sW,\sR)$ 
and their Coxeter elements $\omega_\Gamma$, where $\s \equiv 1$,  have the following special properties
(see, for example, \cite{Bourbaki:Lie} \cite{ACampo:Coxeter} \cite{McMullen:Coxeter}).
\begin{description}
\item{(i)} ({\it Presentation.}) The Coxeter group $\sW$ has presentation in terms of the standard generators $\sR=\{s_1,\dots,s_k\}$:
$$
\langle s_1,\dots,s_k \ : \ (s_is_j)^{m_{i,j}} \rangle.
$$
\item{(ii)} ({\it Monotonicity.}) The spectral radius of Coxeter elements is monotone increasing with respect to inclusion of Coxeter graphs that respect orderings on vertices. If $\Gamma_1$ and $\Gamma_2$ are both higher rank and $\Gamma_1 \subsetneq \Gamma_2$,
then the spectral radius of $\Gamma_2$ is strictly higher than that of $\Gamma_1$.
\item{(iii)} ({\it Bipartite Coxeter eigenvalue.})  The spectral radius of Coxeter elements is bounded from below  
by 
$$
\frac{(\mu^2 - 2) + \sqrt{(\mu^2-2)^2 - 4}}{2},
$$
known as the {\it bipartite eigenvalue} of the Coxeter system.
\end{description}

These properties do not necessarily hold for mixed-sign Coxeter systems (see \cite{Armstrong:thesis}).   

\begin{example}{\em If $v_i$ and $v_j$ are two adjacent vertices on $\Gamma$, and $\s(i) \neq \s(j)$, then $s_is_j$ has infinite order.}
\end{example}

\begin{example} {\em The mixed-sign Coxeter graph $(\Gamma,\s)$, where $\Gamma$ is the complete graph on $3$ vertices,
and $\s \equiv -1$, is spherical.  The Coxeter group is isomorphic to the symmetric group on 4 letters, while the
three generator group with only the pairwise relations  is isomorphic to the affine group associated
to $\widetilde {A_2} = (\Gamma,+1)$ and has infinite order.}
\end{example}

\begin{example}{\em
Property (ii), the monotonicity property, of classical Coxeter elements implies that there is a lower bound greater than 1 for the spectral radius of classical
Coxeter elements  of non-spherical or affine Coxeter graphs.
For classical Coxeter elements with $\mu^2 > 2$, the smallest positive spectral radius is Lehmer's number
$$
\lambda_L \approx 1.17628
$$
(see \cite{McMullen:Coxeter}),
which is also the smallest possible bipartite eigenvalue of a non-spherical and non-affine Coxeter graph. 
On the other hand, as we see later in this paper, the spectral radius of mixed-sign Coxeter elements can be made arbitrarily close to one.
}
\end{example}


\subsection{Representations of Artin groups}\label{Artin-sec}

We recall the definition of the Artin group $\sA_\Gamma$ associated to a Coxeter graph $\Gamma$, and define representations of
$\sA_\Gamma$ associated to an ordered mixed-sign Coxeter graph $(\Gamma,\s)$ (see, for example, \cite{Birman:Braids}
and references therein).

The {\it Artin group} of $\Gamma$ is the group
$$
\sA_\Gamma = \langle \sigma_1,\dots,\sigma_n \ : \ [\sigma_i \sigma_j]_{m_{i,j}}  = [\sigma_j\sigma_i]_{m_{i,j}}\rangle
$$
where $[\sigma_i \sigma_j]_m$ is the alternating product 
$$
[\sigma_i \sigma_j]_m = \sigma_i \sigma_j \sigma_i \dots
$$
of length $m$.
If $\Gamma$ is the classical Coxeter graph $A_n$, then $\sA_\Gamma$ is
the  braid group on the disk with $n+1$ punctures.

Let
$\R^{\sV}$ be the vector space of real labels on the vertices of $\Gamma$.    The ordered vertices $\sV$ determine an ordered basis  
$v_1,\dots,v_n$ of $\R^{\sV}$.
Let $A$ be the adjacency matrix for $\Gamma$.  Let $F$ be the skew-symmetric bilinear
form on $\R^{\sV}$ defined with respect to $v_1,\dots,v_n$ by 
$$
F=A^+ - A^-, 
$$
where 
$A^+$ is the upper triangular part of $A$, and $A^-$ is the lower triangular part.
The matrix $F$ defines a skew symmetric form on $\R^{\sV}$ that depends on the choice of ordering of $\sV$.

Let $\s$ be a sign-labeling for $\Gamma$.
Define $\rho_{\Gamma} $ to be the representation
$$
\rho_{\Gamma}: \sA_\Gamma \rightarrow \GL(\R^{\sV})
$$
 preserving $F$ defined
by 
\begin{eqnarray*}
\rho_{\Gamma}(\sigma_i)(v_j) &= &
v_j + F(v_i,v_j) v_i\\
&=&
\left \{
\begin{array}{rl}
v_j &\qquad \mbox{if $i=j$,}\\
v_j + a_{i,j}  v_i &\qquad \mbox{if $i < j$,}\\
v_j - a_{i,j}  v_i &\qquad \mbox{if $i  > j$.}
\end{array}
\right .
\end{eqnarray*}
Then the image preserves $F$.
We call $\rho_{\Gamma}$ the {\it Artin representation} of $\sA_{\Gamma}$.
Define
$$
\sigma_{\Gamma,\s} = \sigma_1^{\s(1)} \cdots \sigma_k^{\s(k)}
$$
to be the {\it Artin element} associated to the ordered mixed-sign Coxeter graph.

\begin{proposition}\label{ArtinCoxeter-prop}  The Coxeter element $\omega_{\Gamma,\s}$ and
the Artin element $\sigma_{\Gamma,\s}$ are related by
$$
\omega_{\Gamma,\sign} = - \rho_{\Gamma}(\sigma_{\Gamma,\sign}).
$$
\end{proposition}

Let $B_c = U + c U^T$, for $c \in \C \setminus 0$.  Then we have $B = B_1$ and $F = B_{-1}$.  
Now define elements $f_1^{(c)},\dots,f_n^{(c)} \in \GL(\R^{\sV})$ by
\begin{eqnarray*}
f^{(c)}_i(v_j) &=& v_j - \s(i) B_c(i,j) v_i\\
&=& 
\left \{
\begin{array}{lr}
-c v_i &\qquad \mbox{if $i=j$}\\
v_j + \s(i)a_{i,j} v_i &\qquad \mbox{if $i < j$}\\
v_j + c\s(i) a_{i,j} v_i & \qquad \mbox{if $i > j$}\\
\end{array} 
\right .
\end{eqnarray*}
Then $f^{(1)}_i = s_i$ and $f^{(-1)}_i=\rho_{\Gamma}(\sigma_i^{\s(i)})$.

Proposition~\ref{ArtinCoxeter-prop} follows from the following generalization of a result of Howlett
 \cite{Howlett:Coxeter}.

\begin{lemma}\label{Howlett-lem}  Using the above notation
$$
f^{c}_1 \cdots f^{c}_n = -c U^{-1} U^T.
$$
\end{lemma}

We present the generalized proof here.

{\bf Proof of Lemma~\ref{Howlett-lem}.}  
First we notice that
\begin{eqnarray*}
U  f^{(c)}_1 &=&
\left [
\begin{array}{ccccc}
\s(1) & - a_{1,2} &-  a_{1,3}& \dots &- a_{1,n}\\
0 &\s(1) &-a_{2,3}& \dots &-a_{2,n}\\
\dots &&&\\
0 & &&\dots & \s(n)
\end{array}
\right ]
\left [
\begin{array}{ccccc}
-c & \s(1)a_{1,2} && \dots &\s(1) a_{1,n}\\
0 & 1 & 0 &\dots & 0\\
\dots &&&&\\
0 &&&&1
\end{array}
\right ]\\
&=&
\left [
\begin{array}{ccccc}
- c \s(1) &0 && \dots &0\\
0 &\s(1) &-a_{2,3}& \dots &-a_{2,n}\\
\dots &&&\\
0 & &&\dots & \s(n)
\end{array}
\right ]
\end{eqnarray*}
Assume that
\begin{eqnarray*}
U f^{(c)}_1 \cdots f^{(c)}_k &=&
\left [
\begin{array}{cc ccc ccc}
- c \s(1) &0 &&& \dots &&&0\\
-c a_{1,2} &-c \s(1) &0&&&& \dots &0\\
 &\dots&&&&&&\\
-c a_{k,1}&\dots&-c a_{k,k-1} &-c\s({k})&0&&\dots &0\\
0 &&\dots&0&\s(k+1)&-a_{k+1,k+2}&\dots&-a_{k+1,n}\\
&\dots&&&&&&\\
0&&&\dots &&0&\s(n-1)&-a_{n-1,n}\\
0 & &&\dots &&&0& \s(n)\\
\end{array}
\right ] \\
&=&
\left [
\begin{array}{cc}
L_k & 0 \\
0 & U_k
\end{array}
\right ].
\end{eqnarray*}
Multiplying on the right by 
$$
f^{(c)}_{k+1}
=
\left [
\begin{array}{c|c|c}
I & 0 & 0\\

\hline
-\s({k+1})c  a_{k+1,1} \dots 
-\s({k+1})c a_{k+1,k} & -c & \s({k+1})a_{k+1,k+2}
 \dots  \s({k+1}) a_{k+1,n}\\
\hline
0&0&I
\end{array}
\right ]
$$
amounts to replacing the 
$k+1$st row of $U f^{(c)}_1\cdots f^{(c)}_k$ by 
$$
[a_{k+1,1}, \dots, a_{k+1,k}, -\s({k+1}),  0 , \dots, 0].
$$
Thus, 
$$
U f^{(c)}_1 \cdots f^{(c)}_{k+1} = 
\left [
\begin{array}{cc}
L_{k+1} & 0 \\
0 & U_{k+1}
\end{array}
\right ]
$$
and
 $L_n = -c U^T$ as desired.
\qed

For each $c$, the transformation $f = f_1^{(c)} \cdots f_n^{(c)}$ satisfies
$$
f^T B_c f = B_c.
$$

The following consequence of Lemma~\ref{Howlett-lem}  completes the proof of 
Proposition~\ref{ArtinCoxeter-prop}.

\begin{proposition}\label{Howlett-prop}
The mixed-sign Coxeter element and representation of the corresponding
element of the Artin group is given by
$$
\omega_{\Gamma,\s} = -U^{-1} U^T,
$$
and
$$
\rho_{\Gamma}(\sigma_{\Gamma,\s}) = U^{-1}U^T.
$$
\end{proposition}

\subsection{Geometric realization and mixed-sign Coxeter mapping class}\label{geo-sec}

In this section, we define a compact oriented surface $S_{\Gamma}$ from an ordered Coxeter fatgraph $\Gamma$, and
a mapping class $\phi_{\Gamma,\s}$ from $\Gamma$ and a sign-labeling $s$.

Let $\Gamma$ be a Coxeter graph with vertices $\sV$.  A {\it fatgraph} (or {\it ribbon}) structure on $\Gamma$
is a choice of cyclic ordering on the edges emanating from each vertex of $\sV$.  For any embedding of a graph 
$\Gamma$ on a surface $S$ there is a corresponding fat graph structure on $\Gamma$.  Conversely, for any
fat graph structure on $\Gamma$ there is a unique oriented closed surface on which $\Gamma$ embeds as
a fat graph so that the complementary components are disks.

Let
$$
\oo: \sV \rightarrow \{1,\dots,n\}
$$
be the bijection corresponding to the ordering on $\sV$. 
For each $v \in \sV$, let $\sV_v \subset \sV$ be the link of
$v$, i.e., the set of vertices connected to $v$ by an edge.
 Then the fat graph structure of $\Gamma$ is equivalent to a
choice of cycle 
$$
\sigma_v = ({i_1},\dots,{i_k}),
$$
 in the symmetric group $\sS_n$ on $n$ elements for each element $v \in \sV$ with $\sV_v= \{v_{i_1},\dots,v_{i_k}\}$.

Construct a  system of annuli $T_{v}$, for $v \in \sV$, with oriented core curve $\gamma_v$
  so that
\begin{enumerate}
\item $T_v$ and $T_w$ are glued together along a square patch if and only if $v$ and $w$ are connected by an edge;
\item if $\oo(v) < \oo (w)$, then $i_{\mbox{alg}}(\gamma_v,\gamma_w) > 0$; and
\item If $v$ is a vertex,  then for each of the  the core curves $\gamma_w$, for $w \in \sV_v$,
intersect $\gamma_v$ in a cyclic ordering that respects the orientation of $\ell_v$ and cycle $\sigma_v$.
\end{enumerate}
\begin{figure}[htbp] 
   \centering
   \includegraphics[width=1in]{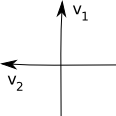}
   \caption{Positive intersection}
   \label{posint-fig}
\end{figure}

Here we use the convention that if $\gamma_i$ and $\gamma_j$ intersect as in Figure~\ref{posint-fig}, then
$i_{\mbox{alg}}(\gamma_i,\gamma_j) = 1$.

\begin{figure}[htbp] 
   \centering
   \includegraphics[width=4in]{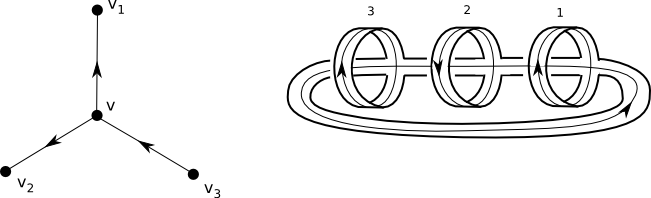}
   \caption{A local picture of Coxeter graph and corresponding union of annuli}
   \label{graphsurface-fig}
\end{figure}

Figure~\ref{graphsurface-fig} shows the arrangement of $T_{v_1},T_{v_2},T_{v_3}$ and $T_{v}$
where  $\sigma_v = (1,2,3)$, and  $\oo(v_2) < \oo(v) < \oo(v_1), \oo(v_3)$.

The arrows in the figure indicate which vertex comes before the other in the global ordering.   One sees that the surface
 depends on the relative global ordering of adjacent vertices and the fatgraph structure.

Let $\s$ be a sign-labeling for $\Gamma$.   Let $\phi_{\Gamma,\s} S_{\Gamma} \rightarrow S_{\Gamma}$ be the
mapping class defined by
$$
\phi_{\Gamma,\s} = (\delta_1)^{\s(1)} \cdots (\delta_k)^{\s(k)},
$$
where $\delta_i$ are the right Dehn twists centered at $\gamma_i$. (See, for example, \cite{FM:MCG} for definition of Dehn twist.)
Then $(S_{\Gamma},\phi_{\Gamma,\s})$ is the {\it Coxeter mapping class} associated to $(\Gamma,\s)$.

\begin{figure}[htbp] 
   \centering
   \includegraphics[width=4in]{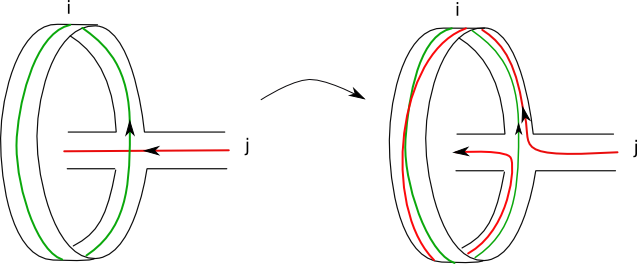}
   \caption{Right Dehn twist along $\delta_i$ acting on $\gamma_j$, where $i< j$.}
   \label{Dehntwist-fig}
\end{figure}

Figure~\ref{Dehntwist-fig} gives an illustration of the action of a right Dehn twist on a transversally intersecting curve.
One may verify that $\delta_i(\gamma_j) = \gamma_j + \gamma_i$ when $\gamma_i$ and $\gamma_j$ are oriented
in this way with $i_{\mbox{alg}}(\gamma_i,\gamma_j) = 1$.

\begin{remark} {\em The surfaces $S_\Gamma$ depend  only on the choice of vertex ordering and fatgraph structure
on $\Gamma$.  Different vertex orderings and fatgraph structures can give rise to different homeomorphism types of surfaces.}
\end{remark}

By contrast, the topological Euler characteristic of $S_{\Gamma}$ is independent of the ordering and fatgraph structure.

\begin{lemma}  Let $S_\Gamma$ be a geometric realization of $\Gamma$.  Then
 $\chi(S_{\Gamma}) = - |\sE|$, 
where $|\sE|$ is the total number of edges of $\Gamma$.  In particular, $\chi(S_\Gamma)$
does not depend on the ordering of $\sV$ or the fatgraph structure.
\end{lemma}

\bold{Proof.} The construction of $S_{\Gamma}$ is inductive with respect to an ordering on the vertices $\sV=\{v_1,\dots,v_k\}$.  
At each stage  $i$ we attach an annulus to the
preceding surface along patches one for each vertex  $v_j$ adjacent to $v_i$  such that $j < i$.  The Euler characteristic
thus changes by the number of such vertices.  Thus each edge is counted exactly once.
\qed

\begin{remark}\label{integer-rem}{\em  We could extend the definition of $(S_\Gamma,\phi_{\Gamma,\sign})$ further, by replacing $\sign$ with 
$$
\epsilon : \sV \rightarrow \Z^*,
$$
where $\Z^*$ is the set of nonzero integers, and setting
$$
\phi_{\Gamma,\epsilon} = \delta_1^{\epsilon_1} \circ \cdots \circ \delta_n^{\epsilon_n}.
$$
Each $\phi_{\Gamma,\epsilon}$ is, however, also associated to a signed graph.  The signed graph $(\Gamma',\sign)$ is obtained
from $(\Gamma,\epsilon)$ by successively replacing each vertex $v_i$ in $\Gamma$ with $m_i=|\epsilon_n|$ copies
$v_i^{(j)}$.  Each of the new vertices $v_i^{(j)}$ of $\Gamma'$ has edges connecting it to all the vertices to 
which $v_i$ was connected.   
The new sign labels on $v_i^{(j)}$ equal the sign of $\epsilon$.
}\end{remark}

Let $\R^{\sV} \hookrightarrow H_1(S_\Gamma;\R)$ be the linear map defined by sending the $i$th basis
vector to $[\gamma_i]$.  Let $GL(\R^{\sV}) \rightarrow GL(H_1(S_\Gamma;\R))$ be the homomorphism defined by extending
by the identity on the complementary space of the image of $\R^{\sV}$.   

\begin{proposition}\label{hom-prop} The induced map $(\phi_{\Gamma,\s})_*: H_1(S_\Gamma;\R) \rightarrow H_1(S_\Gamma;\R)$ on homology
satisfies
$$
(\phi_{\Gamma,\s})_* = \rho_{\Gamma}(\sigma_{\Gamma,\s}) 
$$
and hence 
$$
\lambda_{\mbox{hom}}(\phi_{\Gamma,\s}) = |\omega_{\Gamma,\s}|.
$$
\end{proposition}

\bold{Proof.} Let $g_i = [\gamma_i]$ be the homology classes. 
The choice of orientations and algebraic intersections of $g_1,\dots,g_k$ satisfy
$$
(\delta_i)_*(g_j^{\s(i)})=
\left \{
\begin{array}{ll}
g_j &\qquad\mbox{if $i=j$}\\
g_j + \s(i) g_i &\qquad \mbox{if $i<j$}\\
g_j - \s(i) g_i &\qquad \mbox{if $i>j$.}
 \end{array}
 \right .
 $$
 This can be checked by examining Figure~\ref{Dehntwist-fig}.
Thus, $(\delta_i)^{\s(i)}_*$ restricted to the image of $\R^{\sV}$ in $H_1(S_\Gamma;\R)$ equals
$\rho_{\Gamma}(\sigma_i)$.  \qed

\begin{proposition}\label{pAcrit-prop}  
 If $\Gamma$ is connected and the spectral radius  of 
the Coxeter element $|\omega_{\Gamma,\s}|$
is greater than one, then $(S_{\Gamma,\s},\phi_{\Gamma,\s})$ is pseudo-Anosov, and the dilatation satisfies
$$
\lambda(\phi_{\Gamma,\s}) \ge |\omega_{\Gamma,\s}|.
$$
\end{proposition}

\bold{Proof.}
By the Nielsen-Thurston classification, any mapping class is either periodic, reducible or pseudo-Anosov \cite{Thurston88}.
Since $\Gamma$ is connected, $\omega_\Gamma$ is irreducible, and hence so is the homological action
of $\phi_\Gamma$.  This implies that $\phi_\Gamma$ is not reducible. 

Since $\lambda_{\mbox{hom}}(\phi_{\Gamma,\sign}) > 1$, $\phi_{\Gamma,\sign}$ is not periodic, so it
must be pseudo-Anosov.  The rest follows from the following well-known inequality (see, e.g. \cite{Rykken99}).
$$
\lambda_{\mbox{hom}} (\phi) \leq \lambda_{\mbox{geo}}(\phi).
$$
\qed

\subsection{Mapping classes on closures and interiors}\label{closure-sec}

Let $S$ be a compact surface, we have defined the mapping class group $\Mod(S)$ to be the group of isotopy classes of orientation
preserving self-homeomorphisms of $S$ up to isotopy relative to the boundary of $S$.
Consider the interior $\mbox{int}(S)$ of $S$.   Topologically, this is homeomorphic to a surface with punctures one for each boundary component
of $S$.
Let  $\Mod(\mbox{int}(S))$ be the group of orientation preserving
self-homeomorphisms on $\mbox{int}(S)$ modulo isotopy.  Then the map
\begin{eqnarray}\label{interior-map}
\alpha: \Mod(S) \rightarrow \Mod(\mbox{int}(S))
\end{eqnarray}
defined by restriction
has kernel generated by Dehn twists centered at boundary parallel simple closed curves.

The following results are well-known, and are contained for example in \cite{Birman:Braids}.

\begin{lemma} If $(S,\phi)$ is a pseudo-Anosov element of $\Mod(S)$, then $\alpha(S,\phi)$ is also
pseudo-Anosov, and the dilatations are the same.
\end{lemma}

\bold{Proof.}  If $(S,\phi)$ is pseudo-Anosov, and $(\sF^\pm,\nu^\pm)$ are its associated stable and unstable foliations, then
$(\sF^\pm,\nu^\pm)$ also define stable and unstable foliations for $\alpha(S,\phi)$, and the stretching factor $\lambda$ is
also preserved.
\qed

Let $\overline S$ be the closed surface obtained by filling in the boundary components 
of $S$ with disks.  Then there is a homomorphism
\begin{eqnarray}\label{closure-map}
\beta:  \Mod(S) \rightarrow \Mod(\overline S)
\end{eqnarray}
defined by extending over disks.  This map is neither one-to-one nor onto.  Furthermore, the image of a pseudo-Anosov
mapping class is not necessary pseudo-Anosov.    We will write $(\overline S,\overline \phi) = \beta(S,\phi)$.

\begin{lemma}\label{closure1-lem} If $(S,\phi)$ is pseudo-Anosov, and none of the boundary components are 1-pronged, then 
$(\overline S,\overline \phi)$ is also pseudo-Anosov with dilatation $\lambda(\overline\phi) = \lambda(\phi)$.
\end{lemma}

The idea of the proof is that if there are no 1-pronged boundary components, then the stable and unstable foliations of $(S,\phi)$ 
determine invariant transverse measured foliations for $\overline \phi$ with expansion factor $\lambda^{\pm 1}$ for $\lambda = \lambda(\phi)$
 (see, for example, \cite{HK:braidbounds}, Lemma 2.5).

We can also refine Lemma~\ref{closure1-lem} as follows.
Let $(S,\phi)$ be a mapping class, where $S$ is compact.  For any boundary component $b$ of $S$, let $b=b_1,\dots,b_s$
be the orbit of $b$ under the action of $\phi$.    Let $cl(S,b)$ be the surface obtained by filling in the boundary components
$b_1,\dots,b_s$ with disks, and let $cl(\phi,b)$ be the extension of $\phi$.   Let $cl(S,\phi,b) = (cl(S,b),cl(\phi,b))$.   If $(S,\phi)$
is pseudo-Anosov, and $b$ is $m$-pronged, then all orbits of $b$ are also $m$-pronged.

\begin{lemma}\label{closure-lem}  If $(S,\phi)$ is pseudo-Anosov, then $(cl(S,b),cl(\phi,b))$ is pseudo-Anosov if
the boundary components $b_1,\dots,b_s$ are not 1-pronged.  In this case, 
$$
\lambda(cl(\phi)) = \lambda(\phi).
$$
\end{lemma}

The proof is the same as for Lemma~\ref{closure1-lem}.

\subsection{Bipartite graphs}\label{bipartite-sec}

In this section we collect some special properties of mixed-sign Coxeter systems, and mixed-sign Coxeter mapping
classes associated to bipartite graphs.

A Coxeter graph $\Gamma$ is {\it bipartite (with bipartite ordering)} if the following hold:
\begin{description}
\item {(i)} its vertices can be separated into two disjoint sets $\sV = \sV_1 \cup \sV_2$ where
the subgraph of $\Gamma$ generated by $\sV_i$ has no edges for $i=1,2$; and
\item{(ii)} by the ordering on $\sV$ has the property that the elements of $\sV_1$ proceed all the elements
of $\sV_2$.
\end{description}

A graph $\Gamma$ is bipartite if and only if
it contains no odd cycles.
Given a sign-labeling $\s$ of a Coxeter graph $\Gamma$.  Let $\overline \s$ be the sign-labeling 
defined by $\overline \s(v) = - \s(v)$ for all $v \in \sV$.  

\begin{theorem}\label{bipartite-thm} If $\Gamma$ is bipartite, and $\s$ is any sign-labeling on $\Gamma$, then
$(\sW_{\Gamma,\s},\sR_{\Gamma,\s})$ and  $(\sW_{\Gamma,\overline\s},\sR_{\Gamma,\overline \s})$ 
are conjugate as subgroups of $GL(\sR^{\sV})$, and, in particular, the spectral radius of Coxeter elements satisfies
$$
|\omega_{\Gamma,\s}| = |\omega_{\Gamma,\overline \s}|.
$$
\end{theorem}

\bold{Proof.}  It suffices to show that the generating sets
$\sR_{\Gamma,\s} = \{s_1,\dots,s_k\}$ and $\sR_{\Gamma,\overline\s} = \{\overline s_1,\dots,\overline s_k\}$
are conjugate as elements of $GL(\R^{\sV})$.  

 let $\sV = \sV_1 \cup \sV_2$ be the bipartite partition.  Let $k_i$ be the number elements
in $\sV_i$, for $i=1,2$, and let $k=k_1  + k_2$ be the total number of vertices $\sV$.  
Let $I_{k_1,k_2}$ be the $k \times k$ diagonal matrix with the first $k_1$ diagonal entries equal to 1
and the second $k_2$ diagonal entries equal to -1.
Then $s_i$ and $s_i'$ satisfy
$$
s_i  = I_{k_1,k_2} \overline s_i I_{k_1,k_2}.
$$
\qed

\begin{figure}[htbp]
\begin{center}
\includegraphics[height=1.25in]{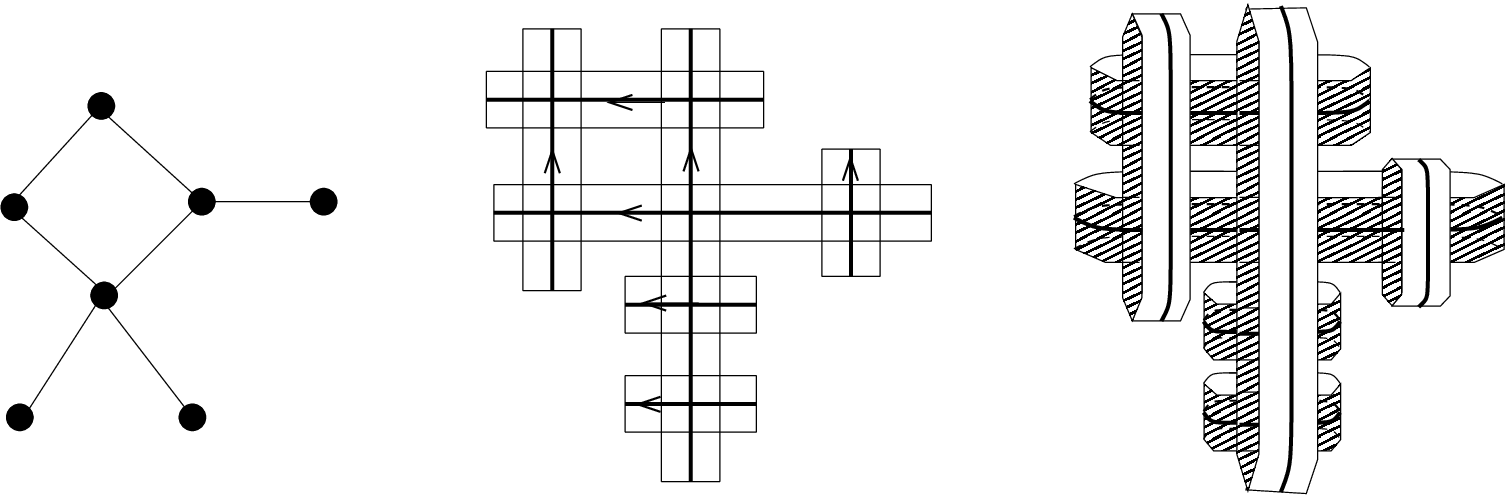}
\caption{Surface associated to a bipartite graph}
\label{bipartite-fig}
\end{center}
\end{figure}

 If a graph $\Gamma$ is bipartite, and is given the bipartite ordering, then there is a fatgraph structure
on $\Gamma$ so that after cutting each annulus at a transversal arc between the two boundary components,
the surface can be placed on a plane as a union of rectangles oriented in vertical and horizontal directions
as in Figure~\ref{bipartite-fig}.  The right diagram in Figure~\ref{bipartite-fig} gives the corresponding 
surface $S_\Gamma$.

The {\it bipartite eigenvalue} of a graph $\Gamma$ is defined by
$$
\beta_\Gamma = |x^2 - (2-\mu^2)x + 1|,
$$
where $\mu$ is the spectral radius of the adjacency matrix of $\Gamma$.
Note, this does not depend on the ordering of the vertices of $\Gamma$.

The following theorem was proved  for (positively signed) classical Coxeter graphs $\Gamma$
in \cite{McMullen:Coxeter}.
\begin{figure}[htbp]
\begin{center}
\includegraphics[height=1in]{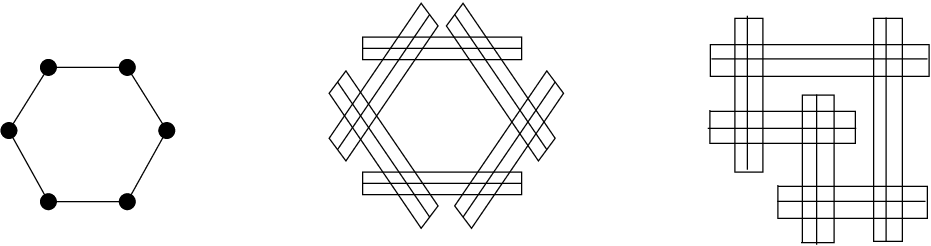}
\caption{Two surfaces obtained from the hexagonal graph by identifying opposite
short edges of the rectangles.}
\label{hexagon-fig}
\end{center}
\end{figure}

\begin{theorem}\label{orientable-thm} Let $\Gamma$ be a positively signed bipartite Coxeter graph
with bipartite ordering.
Let $\mu$ be the spectral radius of the adjacency matrix of $\Gamma$.   Then either $|\mu^2-2| \leq 2$, which
implies that $|\omega_{\Gamma,1}| = 1$, or 
$$
|\omega_{\Gamma,1}| = \beta_\Gamma
$$
and hence only depends on the combinatorics of $\Gamma$.
Furthermore, for any arbitrary (positively signed) Coxeter graph $\Gamma$, with
arbitrary ordering on the vertices,
$$
 |\omega_{\Gamma,1}| \ge \beta_{\Gamma}.
$$
\end{theorem}

In \cite{Thurston88}, Thurston gave an example of pseudo-Anosov mapping classes
constructed using classical bipartite Coxeter graphs.   These are mixed-sign Coxeter
mapping classes associated to a bipartite Coxeter graph, with bipartite order.  He proved
the following.

\begin{theorem}  If $\Gamma$ is a classical bipartite Coxeter graph with bipartite order, then 
the widths and lengths of the rectangles in the construction of $S_{\Gamma}$ can be
chosen so that 
$\phi_{\Gamma,\s}$  has constant derivative.  If $\phi_{\Gamma}$ is pseudo-Anosov, then 
$$
\lambda_{\mbox{geo}}(\phi_{\Gamma,\s}) = \lambda_{\mbox{hom}}(\phi_{\Gamma,\s}) = \beta_{\Gamma}.
$$
\end{theorem}

\begin{figure}[htbp]
\begin{center}
\includegraphics[height=0.75in]{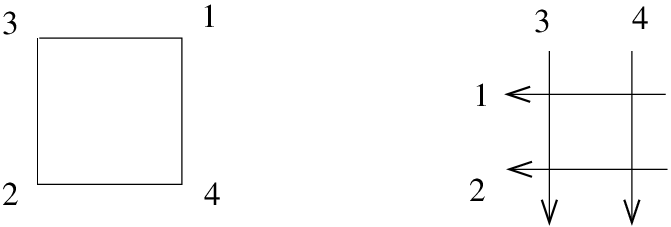}
\caption{Dual configuration associated to a 4 cycle with bipartite ordering}
\label{bsquare-fig}
\end{center}
\end{figure}

\begin{example}{\em Consider the hexagonal graph shown $\Gamma_6$
in Figure~\ref{hexagon-fig}.  Since
each vertex has order two, there is only one fatgraph structure on $\Gamma_6$.
For the middle
diagram, there is no way to orient the core curves on the rectangles in a way
that is orientation compatible with any ordering on $\Gamma_6$.   The orientation
on one curve determines the orientations on its adjacent ones.   Thus,
a $2n$-gon has an orientation compatible diagram of this form if and only if $n$
is even.  The right 
diagram is compatible with the bipartite ordering.  
In this example, the middle surface has genus
2 and 4 boundary components, while the right hand surface has genus 3 and 2 boundary 
components.
}\end{example}

\begin{figure}[htbp]
\begin{center}
\includegraphics[height=0.75in]{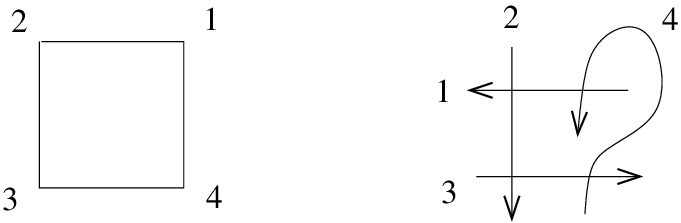}
\caption{Dual configuration associated to a 4 cycle with cyclic ordering}
\label{square-fig}
\end{center}
\end{figure}

\begin{example}{\em It also may not be possible to make the surface from
a globally planar configuration of straight paths, as one can in the bipartite case.  Consider for example the cyclic graph with
4 vertices.  The bipartite ordering gives rise to a planar configuration (see Figure~\ref{bsquare-fig}),
while for the cyclic ordering, one can verify that there is no   configuration of straight line paths that realize the graph
(see Figure~\ref{square-fig}).  This example also illustrates that while the graph determines
the Euler characteristic of the surface (with boundary), but ordering of the vertices can affect
the genus.  In the bipartite case, the surface has $(g,n) =(1,4)$, while in the cyclic case
the surface has $(g,n) = (2,2)$.
}\end{example}


\subsection{Minimum dilatation orientable examples.}\label{orientable-sec}

Table~\ref{graph-table} displays mixed-sign Coxeter graphs that give rise to the minimum dilatation orientable mapping 
classes of genus $2$ through $5$ found by Lanneau and Thiffeault \cite{LT09}.   For genus $2$ and $3$, the equivalent  integer labeled graphs are also given (see Remark~\ref{integer-rem}).
We use the convention that a filled in vertex is given the sign label `+1' while the unfilled vertex is given the sign
label `-1'.

\begin{table}\label{smallorientable-table}
\begin{center}
{\renewcommand{\arraystretch}{2} \renewcommand{\tabcolsep}{0.5cm}
\begin{tabular}{ m{0.4in} | m{4in} }

genus & mixed-sign Coxeter graph\\
\hline
&\\
2 &
\includegraphics[width=2in]{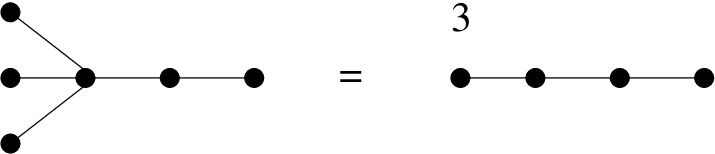}
\\
3
&
\includegraphics[width=2.75in]{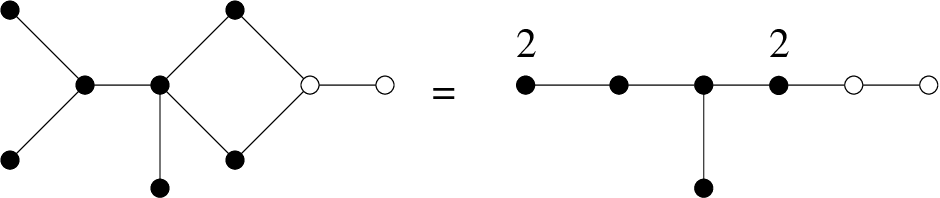}
\\
4
&
\includegraphics[width=1.75in]{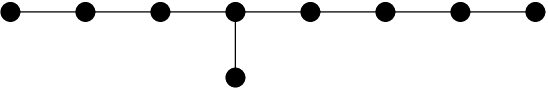}\\

5
&
\includegraphics[width= 2in]{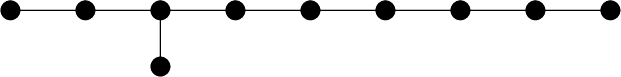}\\

\end{tabular}
}
\caption{Graphs corresponding to minimum dilatation orientable examples for $g=2,3,4,5$}\label{graph-table}
\end{center}
\end{table}

Let $\delta_g^+$  be the minimum dilatation for an orientable pseudo-Anosov mapping class
on a closed surface.
The minimum orientable examples for genus $g=2,4,5$ can be realized as the closures of classical Coxeter mapping
classes associated to bipartite graphs constructed in  \cite{Thurston88} \cite{Leininger04}.  
For genus 4 and genus 5, the graphs are the classical hyperbolic extensions of the $E_7$ and $E_8$ graphs,
and in the genus 5 example, the dilatation is equal to Lehmer's number.
The hyperbolic extension of $E_6$ is given in Figure~\ref{genus3alt-fig}.   It's corresponding mapping class
is defined on a surface of genus $4$, but it has the same dilatation as the genus 3 example given in Table~\ref{graph-table}.

The sequence $\delta_g^+$ converges to 1.  Moreover, we have
$$
\log (\delta_g^+) \asymp \frac{1}{g}
$$
(see \cite{HK:braidbounds}).  
Thus, classical Coxeter mapping classes cannot realize small dilatation orientable mapping 
classes for high genus.

\begin{remark}\label{genus3-rem}{\em One can also study dilatations of pseudo-Anosov mapping classes $(S,\phi)$ in terms of the topological Euler characteristic 
$\chi(S)$.   This is well-motivated by the following. Given a pseudo-Anosov mapping class $(S,\phi)$, let 
$$
L(S,\phi) = \lambda(\phi)^{|\chi(S)|}
$$
be the {\it $\chi$-normalized dilatation} of $(S,\phi)$.  
Let $M$ be
a hyperbolic 3-manifold, and $F$ a fibered face (see \cite{Thurston:norm}).  Then $F$ is a polyhedron of dimension equal to $b_1(M) - 1$
where $b_1(M)$ is the first Betti number of $M$.  The rational points in the interior of $F$ correspond to mapping classes
$(S,\phi)$ that are monodromies of fibrations of $M$ over the circle. 

On any compact subset of the interior of a fibered face $F$
the $\chi$-normalized dilatation 
extends to a continuous convex function on $F$ (see \cite{Fried82} \cite{McMullen:Poly}) and hence is bounded.
Given a pseudo-Anosov $(S,\phi)$, let $(S^0,\phi^0)$ be the mapping class obtained by letting $S^0 = S \setminus \mbox{Sing}(\phi)$,
and $\phi^0 = \phi |_{S^0}$.  Then $(S^0,\phi^0)$ is pseudo-Anosov, and $\lambda(\phi^0) = \lambda(\phi)$ (see, Lemma~\ref{closure-lem}).

Farb, Leininger and Margalit Universal Finiteness Theorem \cite{FLM09} implies that
the collection of mapping classes $(S,\phi)$ with bounded {\it $\chi$-normalized dilatation}
correspond (after puncturing $S$ at singularities) to rational points on compact subsets of a finite collection fibered faces.   
 \begin{figure}[h]
\begin{center}
\includegraphics[width=1.4in]{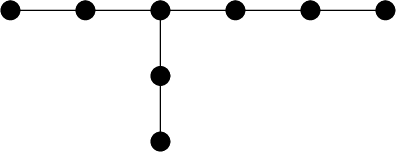}
\caption{A positive Coxeter graph related to genus 3 example.}
\label{genus3alt-fig}
\end{center}
\end{figure}

As an example, numerically
$\dil_3^+$  equals the house of the hyperbolic extension of $E_6$ (see \cite{McMullen:Coxeter}, Table 5, and Figure~\ref{genus3alt-fig}). 
Let $(S_{hE_6},\phi_{hE_6})$ be the mapping class obtained from $hE_6$ using the bipartite ordering.  Then $S_{hE_6}$ 
has genus 4 and 2 boundary components.   The genus 3 example obtained from the graph $\Gamma_3$ in Table~\ref{graph-table} has 4 
boundary components.
Thus, the surfaces $S_{hE_6}$ and $S_{\Gamma_3}$ have the same topological Euler characteristic.     
}\end{remark}

The above discussion suggests another version of the minimum dilatation problem.

\begin{problem} Find the minimum dilatation of pseudo-Anosov mapping classes with no interior singularities with a given
topological Euler characteristic.
\end{problem}

%

\section{Twist graphs and twisted Coxeter mapping classes}\label{twist-sec}

In this section, we define (full) twist graphs, and associated (full) twist mapping classes and
use them to prove Theorem~\ref{main2-thm}.  Twist graphs are elementary 
building blocks  that can be used to construct small dilatation pseudo-Anosov mapping classes via
Murasugi sum.   We will construct sequences of mapping classes
associated to iterative joins of twist graphs, and investigate conditions under 
which the normalized dilatations are bounded.

Before defining twist graphs and twist mapping classes, we define Murasugi sum for a pair fibered 3-manifolds
(Section~\ref{Murasugi-sec}).  As as example, we recall the definition of Hopf-plumbing for fibered knot and
link complements (Section~\ref{Hopf-sec}).

%
%
%

\subsection{Generalized Murasugi sums of mapping classes.}\label{Murasugi-sec}

The Murasugi sum was originally defined for fibered links in $S^3$ \cite{Murasugisum}.
In this section we study properties of Murasugi sums for arbitrary mapping classes.

Let $P_{2k}$ be a $2k$-sided polygon with alternate edges removed. 
The polygon $P_{2k}$ is {\it properly embedded} in a compact surface $S$ if
the boundary components of $P_{2k}$ are contained in the boundary
of $S$, and the interior of $P_{2k}$ is contained in the interior of $S$.

Let $(S_0,\phi_0)$ and $(S_1,\phi_1)$ be two mapping classes
with proper embeddings of $P_{2k}$.   Let $S$ be the surface obtained
by gluing $S_0$ and $S_1$ by identifying the interiors of the embedding
of $P_{2k}$ in $S_0$ to the interior of the embedding of $P_{2k}$ in
$S_1$ after rotating by $\frac{2\pi}{k}$.

Note that the intersection of the closure of $S \setminus S_i$ with $S_i$ is contained in
the boundary of $S_i$ for $i=0,1$.
Thus, we can 
extend $\phi_i : S_i \rightarrow S_i$ by the identity on $S \setminus S_i$ and let
$\phi : S \rightarrow S$ be the composition $\phi= \phi_1 \circ \phi_0$.
The mapping class $(S,\phi)$ is called the {\it Murasugi sum} of 
$(S_0,\phi_0)$ and $(S_1,\phi_1)$ relative to the embeddings of $P_{2k}$.

\bold{Murasugi sum and mapping tori.}  Let $M_0$ and $M_1$ be the mapping tori
of $(S_0,\phi_0)$ and $(S_1,\phi_1)$.  Identify $S_i$ with a fiber of $M_i$.
Let $M_i^\lozenge$ be the result of cutting $M_i$ along $\iota_i(P_{2k}) \subset S_i$,
creating a boundary component homeomorphic to a sphere with hemispheres identified
with two copies $\iota_i(P_{2k}^{+})$ and $\iota_i(P_{2k}^{-})$  of $P_{2k}$ glued
together along their boundaries.  
Let $M'$ be the result of gluing $M_0^\lozenge$ with $M_1^\lozenge$ along boundary
spheres, so that
 $\iota_0(P_{2k}^+)$  is glued to $\iota_1(P_{2k}^-)$ and
 $\iota_0(P_{2k}^-)$  is glued to $\iota_1(P_{2k}^+)$.

\begin{lemma}\label{Murasugi-lem} The mapping torus $M$ of $(S,\phi)$ is 
homeomorphic to $M'$.
\end{lemma}

\begin{figure}[htbp] 
   \centering
   \includegraphics[height=1.75in]{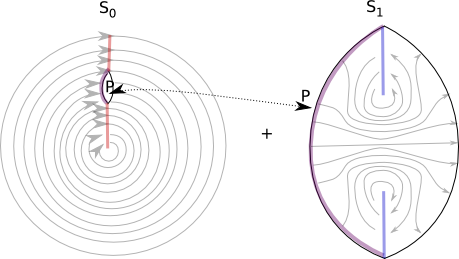} 
   \caption{Murasugi sum of surface flows.}
   \label{mappingtorus-fig}
\end{figure}

To prove Lemma~\ref{Murasugi-lem} it is useful to view the mapping torus of
a mapping class $(S,\phi)$ as a manifold with a 
continuous surjection
$$
\mathfrak f: S \times [0,1] \rightarrow M
$$
 such that
\begin{description}
\item {(i)} $\mathfrak f$ is  an embedding on $S \times [1,0)$, and
\item {(ii)} $\mathfrak f(s,1) = \mathfrak(\phi(s),0)$, for $s \in S$.
\end{description}
An $f$ satisfying (i) and (ii) is called a {\it surface flow}, with {\it transverse surface} $S$
identified with\linebreak $\mathfrak{f} (S \times \{0\})$,
and {\it monodromy} $(S,\phi)$.

\begin{lemma}\label{surfaceflow-lem}  The map from fibrations to surface flows given by
modding out the product $S \times [0,1]$ by the equivalence $(x,1) \sim (\phi(x),0)$
is a bijection.
\end{lemma}

\bold{Proof.} Given a surface flow with transverse surface $S$ and monodromy 
$(S,\phi)$,  there is corresponding
fibration of $M$ over $S^1$ with monodromy $(S,\phi)$ defined by contracting the embedded
surfaces $\mathfrak f (S \times \{t\})$.    \qed

We relax the definition of surface flow slightly  to include the following.
Let $f : S \times [0,1] \rightarrow M$ be a continuous surjective map so
that $(i)$ is replaced by
\begin{description}
\item{(ia)} $\mathfrak f |_{(S \times \{t\})}$ is 1-1 for all $t \in [0,1]$; and
\item{(ib)} $f(s,t) = f(s',t')$ only if $s = s'$ and for all $t_1 \in [t,t']$, $f(s,t_1) = f(s,t)$.
\end{description}

 Given a flow $f$ satisfying $(ia),(ib)$ and $(ii)$, 
one can continuously deform $f$ until it satisfies $(i)$ and $(ii)$, and hence
$f$ determines a unique fibration of $M$ with fiber $S$.  We will use this weaker
version of surface flow in what follows.

\begin{figure}[htbp] 
   \centering
   \includegraphics[height=1.5in]{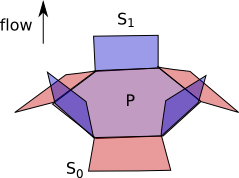} 
   \caption{Local gluing of $S_0$ and $S_1$ at $P$ with upward surface flow.}
   \label{polygon-fig}
\end{figure}

\bold{Proof of Lemma~\ref{Murasugi-lem}.}
Let $P = P_{2k}$, and let 
 $(M_i,S_i,\f_i)$ be two flows defined (up to isotopy) by $(S_i,\phi_i)$, for $i=0,1$..
Define a flow $(M,S,\f)$ by
$\f: S \times [0,1] \rightarrow M$, where
$$
f(s,t) =
\left \{
\begin{array}{ll}
\f_0(s,2t) &\qquad{\mbox{if $s \in S_0, \ 0 \leq t \leq 1/2$}}\\
\f_1(s,0) &\qquad\mbox{if $s \in S_1 \setminus \iota_1(P), \ 0 \leq t \leq 1/2$}\\
\f_0(s,1) &\qquad\mbox{if $s \in S_0 \setminus \iota_0(P), \ 1/2 \leq t \leq 1$}\\
\f_1(s,2t-1) &\qquad\mbox{if $s \in S_1,\ 1/2 \leq t\leq 1$}
\end{array}
\right .
$$

The sum of the flows is illustrated in Figure~\ref{mappingtorus-fig}.
The monodromy $\phi$ is the isotopy type of the composition $f_1 \circ f_0$.
Figure~\ref{polygon-fig} illustrates the local gluing of $S_0$ and $S_1$.
The map $f$ defines a unique surface flow on $M$ up to isotopy, and hence a fibration of 
$$
M  \rightarrow S^1,
$$
with monodromy equal to the Murasugi sum of $(S_0,\phi_0$ and $(S_1,\phi_1)$.
 \qed
 
\begin{figure}[htbp] 
   \centering
   \includegraphics[width=4in]{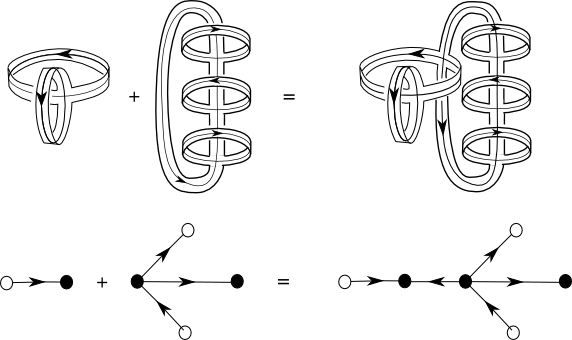} 
   \caption{Join of graphs, and corresponding geometric realization.}
   \label{join-fig}
\end{figure}
 
\begin{example}{\em  Let $\Gamma_i = (\sV_i,\sE_i)$, $i=0,1$ be two mixed-sign Coxeter graphs, and
let $v_i \in \sV_i$ be fixed vertices.    Then the {\it join} of $\Gamma_0$ and $\Gamma_1$ at $v_0$ and $v_1$ is the
graph obtained by taking the disjoint union of $\Gamma_0$ and $\Gamma_1$ and adding a new edge $\epsilon$ between $v_0$ and $v_1$.
If $\Gamma_0$ and $\Gamma_1$ are ordered and have fatgraph structure at their vertices, then the {\it join} is
determined by the choice of element $w_i \in \sV_{v_i}$ for each $i=1,2$.  Then $v_1$
is inserted into $\sV_{v_2}$ after $w_2$, and similarly $v_2$ is inserted into $\sV_{v_1}$
after $w_1$.

Assume that $\Gamma_i$ are simply-laced Coxeter graphs with global orderings and local fatgraph structure, 
 let $\s_i$, $i=0,1$ be sign labels on $\sV_i$, and let $(S_{\Gamma_i,\s_i},\phi_{\Gamma_i,\s_i})$
be the corresponding Coxeter mapping class.  Then for any choice of pairs of adjacent edges at $v_0$ and $v_1$,
we have a new mixed-sign Coxeter graph with global and local fatgraph structure $(\Gamma,\s)$, and the
Coxeter mapping class $(S_{\Gamma,\s},\phi_{\Gamma,\s})$ is obtained from $(S_{\Gamma_i},\s_i)$ by 
Murasugi sum along square a $P_4$.

Figure~\ref{join-fig} shows an example.   Although the graphs are drawn with sign-labelings, the geometric realizations
only depend on the underlying ordered fatgraphs.
}
 \end{example} 
 
 \subsection{Hopf plumbing for fibered links}\label{Hopf-sec}

Since the disk is a fiber surface for the complement of the unknot, any sequence of Hopf-plumbings starting
with the unknot gives rise to a new fibered 3-manifold.    Furthermore, the monodromy is simply the composition of the old monodromy
with the monodromy of the Hopf link, that is a right Dehn twist if the twist is clockwise, and a left Dehn twist if
the twist is counter-clockwise. 

In many cases the mixed-sign Coxeter mapping classes are the monodromy of fibered knots and links in $S^3$. 
For example, if we take a configuration of oriented chords $\ell_1,\dots,\ell_k$ on an oriented disk in $S^3$
such that the algebraic intersections $\iota_{\mbox{alg}}(\ell_i,\ell_j) > 0$ for $i > j$, then we can construct
a surface $S_{\Gamma,\s} \subset S^3$  for any choice of labeling $\s$ as follows.
For $i=1,\dots,k$, we successively attach a band with a coutner-clockwise (resp., clockwise)
full-twist along each chord\ $\ell_i$ according to whether $\s(i)$ is positive (resp., negative).  This construction
was studied in the positive case ($\s \equiv 1$) in \cite{Hironaka:Coxeter}.

\begin{proposition} The link complement
$S^3 \setminus \partial S_{\Gamma,\s}$ fibers over the circle with fiber $S_{\Gamma,\s}$ and monodromy
$(S_\Gamma,\phi_{\Gamma,\s})$.
\end{proposition}

\bold{Proof.} By construction, the link $K = \partial S_{\Gamma,\s}$ is obtained by a sequence of Hopf-plumbings starting with
the disk.   Thus, $S^3 \setminus K$ is fibered with fiber $S_{\Gamma,\s}$.

\begin{figure}[htbp] 
   \centering
   \includegraphics[width=4in]{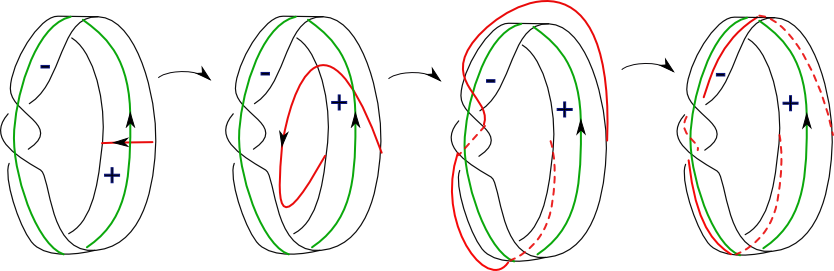}
   \caption{The monodromy of a counter-clockwise twisted Hopf band.}
   \label{Hopftwist-fig}
\end{figure}

To find the monodromy, we need only verify that a clockwise twisted Hopf band has monodromy equal to a right Dehn twist.
This is shown in Figure~\ref{Hopftwist-fig}.
\qed

\bold{Seifert matrices.}
Later in this paper, it will be useful to write down the homological dilatation matrix for the monodromy of
a fibered link complement.    

Let $K \subset S^3$, and let $S$ be a surface embedded in $S^3$ with $K = \partial S$.   Let $x_1,\dots,x_k$ be
a system of simple closed curves on $S$ that generate $H_1(S:\Z)$.   We can find a matrix representing the action of
the monodromy of the fibration as follows (see \cite{Rolfsen76}).   For each $i = 1,\dots,k$, let $x_i^+$ be the
result of pushing $x_i$ into $S^3 \setminus S$ in the positive direction.  
Let $A$ be the $k \times k$ matrix whose $i,j$ entry is given by the linking number $\ell k(x_i^+,x_j)$.   

\begin{proposition}\label{Seifert-matrix-prop}
The linear map $\phi_*: H_1(S;\Q) \rightarrow H_1(S;\Q)$ induced by the monodromy $(S,\phi)$ is defined by the matrix
$$
\phi_*: (A^{\mbox{tr}})^{-1} A.
$$
\end{proposition}

The matrix $(A^{\mbox{tr}})^{-1} A$ is sometimes known as the Seifert matrix for the spanning surface $S$.

\subsection{Fibered faces and dilatations.}\label{fiberedface-sec}

Thurston's fibered face theory for 3-manifolds provides a way to parameterize mapping classes with related dynamics.
Let $M$ be an oriented 3-manifold.  An element of $\alpha \in H^1(M;\Z)$ is {\it fibered} if the corresponding map
$\alpha_*: \pi_1(M) \rightarrow \Z$ has finitely generated kernel.   By a theorem of Stallings \cite{Stallings:fibered}, this
is equivalent to the existence of a fibration $M \rightarrow S^1$ inducing the map $\alpha_*$ on fundamental groups.
If $\alpha$ is fibered, let $(S_\alpha,\phi_\alpha)$ be the {\it monodromy} of $\alpha$, that is, $S_\alpha \subset M$ is a general
fiber, and $\phi_\alpha: S_\alpha \rightarrow S_\alpha$ is the first return map under the flow defined by $\alpha$.

Thurston defined a norm on $H^1(M;\R)$ as follows.   For $\alpha \in H^1(M;\Z)$, let 
$$
||\alpha|| = \min\{|\chi(S_\alpha)|\}
$$
where $S_\alpha$ ranges over  oriented surfaces in $S$
dual to $\alpha$ after removing any connected components of positive Euler characteristic.  Thurston showed the following.

\begin{theorem}[Thurston \cite{Thurston:norm}]   If $M$ is hyperbolic, then $|| \ ||$ extends to a norm on $H^1(M;\R)$ and the unit ball is a compact convex polyhedron.
Furthermore, for every top dimensional open face $F$ of the unit ball, either there are no fibered elements in the cone $F \cdot \R^+$
or all the integral points in $F \cdot \R^+$ are fibered with pseudo-Anosov monodromy.
\end{theorem}

If the cone $F \cdot \R^+$ contains fibered elements, it is called a {\it fibered cone} and $F$ is a {\it fibered face}.  Given a rational
ray from the origin passing through a fibered face $F$, there is a unique {\it primitive} integral element $\alpha$ on the ray with relatively prime
coordinates.  This element $\alpha$ corresponds to a fibration of $M$ over the circle with connected fiber $S_\alpha$, and for any positive integer $k$,
the element $k \alpha$ corresponds to a fibration of $M$ whose fiber is $k$ copies of $S_\alpha$.  Furthermore, the dilatations of
the monodromy $\phi_\alpha$ and $\phi_{k\alpha}$ are related by
$$
\lambda(\phi_{k\alpha})  = \lambda(\phi_\alpha)^{\frac{1}{k}}, 
$$
or
$$
\log \lambda(\phi_{k_\alpha}) = \frac{1}{k} \log \lambda(\phi_\alpha).
$$

\begin{theorem}[Fried \cite{Fried82}] The function 
$$
l (\alpha) = \log \lambda(\phi_\alpha)
$$
extends to a continuous convex function on each fibered cone $F  \cdot \R^+$ that is homogeneous of degree $-1$
and goes to infinity toward the boundary of $F$.
\end{theorem}

\begin{corollary}  For $F$  a fibered face, and $\alpha \in F \cdot \R^+$,  let $\overline\alpha = \frac{\alpha}{||\alpha||}$.
Then the normalized dilatation function
$$
L(\overline\alpha) = \lambda(\phi_\alpha)^{|\chi(S_\alpha)|}
$$
defined for integral elements extends to a continuos and convex function on $F$ that goes to infinity toward the boundary of $F$.
\end{corollary}

\bold{Proof.} This follows from the fact that when $\alpha$ is a fibered element $||\alpha|| = |\chi(S_\alpha)|$,
where $S_\alpha$ is the fiber of the corresponding fibration of $M$ to $S^1$.
\qed

\begin{corollary} For $F$ a fibered face, and $K \subset F$ an infinite compact subset, the collection of monodromies $(S_\alpha,\phi_\alpha)$
corresponding to rational points $\overline \alpha \in K$ has unbounded topological Euler characteristic and bounded normalized dilatation.
In particular, if $\overline \alpha_n$ is a  sequence of distinct rational points on $F$ converging to an interior element in $F$, then
the corresponding monodromies $(S_n,\phi_n)$ have unbounded normalized dilatations and bounded normalized dilatation.
\end{corollary}

\subsection{Full twist braids and their monodromy.}\label{twists-sec}

\begin{figure}[htbp] 
   \centering
   \includegraphics[height=1in]{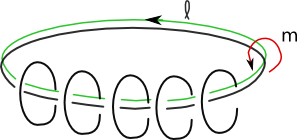}  
   \caption{The mapping torus associated to $(\Sigma_m^{(k)},R_m^{(k)})$.}
   \label{link-fig}
\end{figure}

In this section, we define twist maps and twist graphs, and use them as building blocks for constructing sequences of 
mixed-sign Coxeter mapping 
classes.   The twist maps $(\Sigma_m^{(k)},R_m^{(k)})$, $m\ge2$, $k \ge 1$, have the following properties:
\begin{enumerate}
\item $(\Sigma_m^{(k)},R_m^{(k)})$ is the Coxeter mapping class associated to a graph with $\s \equiv -1$;
\item $(R_m^{(k)})^{km}$ is a product of left Dehn twists on the boundary components of $\Sigma_m$;
\item the mapping class $R_m^{(k)}$ preserves a flat structure on $\Sigma_m^{(k)}$ with 
a distinguished periodic orbit $\sO$ of order $m$; and 
\item the mapping tori of the 
 restrictions $(\Sigma_m^{(k)}\setminus \sO, R_m^{(k)}|_{\Sigma_m^{(k)}\setminus \sO})$ are independent of $k$,
 and are homeomorphic to the complement of a tubular neighborhood of the link in $S^3$ drawn in Figure~\ref{link-fig}.
 \end{enumerate}

\begin{figure}[htbp] 
   \centering
   \includegraphics[height=2in]{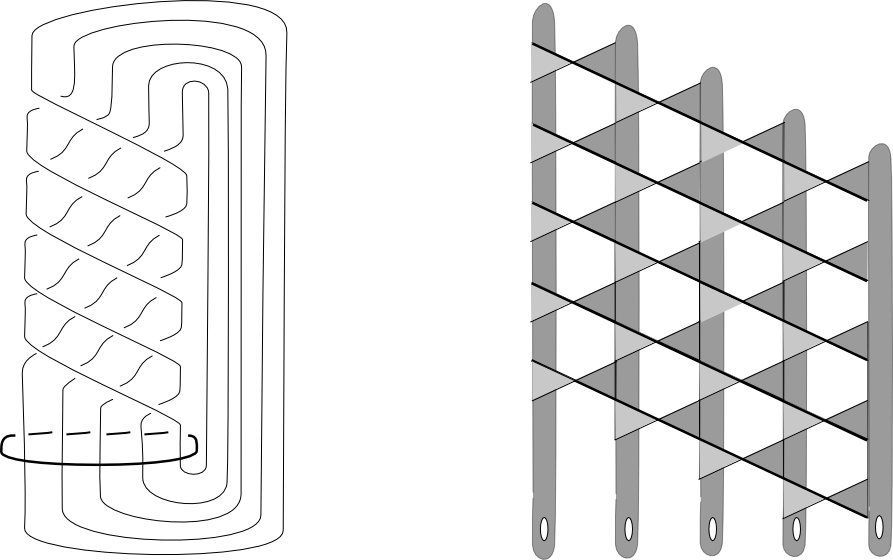}  
   \caption{Closure of a full twist braid on 5 strands, and a fiber surface.}
   \label{twist-fig}
\end{figure}

Let $b_m^{(k)}$ be the product of $k$ full twist braids on $m$ strands.  The corresponding encircled link $L_m^{(k)} \cup E$,
where $E$ is the encircling link,
is drawn in Figure~\ref{twist-fig} in the case
$m=5$ and $k=1$. 

\begin{lemma}  The link complement $S^3 \setminus L_m^{(k)} \cup E$ is independent of $k$.
\end{lemma}

\bold{Proof.}  The complement 
of a tubular neigborhood $L_m^{(k)} \cup E$ in 
$S^3$ is homeomorphic to the complement of the link drawn in Figure~\ref{link-fig}  The link complement for $L_m^{(k)}$
is obtained from the link in Figure~\ref{link-fig} by the surgery that contracts the curve $k \ell + m$,
where $m$ is the meridian and $\ell$ is the longitude of the encircling component $E$.   \qed

Consider the fibration of $S^3 \setminus L_m^{(k)} \cup E$ with fiber equal to the surface $\Sigma_m^{(k)}$ drawn on the right in Figure~\ref{twist-fig}.
The monodromy is periodic of order $km$ (modulo the action near the boundary), and can be seen  explicitly as follows.  

The surface $\Sigma_m^{(k)}$ is a union of $m$ {\it main} disks (drawn vertically in Figure~\ref{twist-fig}) and $km$ {\it attaching} disks (drawn horirzontally).
We can think of the darkly shaded regions as being the positively oriented
side  of $\Sigma_m^{(k)}$, and the lightly shaded region as being on the negatively oriented side.  Number the main disks from right to left
$d_1,\dots,d_m$ and the attaching disks $a_1,\dots,a_{km}$ from top to bottom.  Then the monodromy acts by cyclically permuting the main disks  
$$
d_1 \rightarrow d_2 \rightarrow \cdots \rightarrow d_m \rightarrow d_1
$$
while rotating them by $\frac{2\pi}{km}$,
and the attaching disks
$$
a_1 \rightarrow a_2 \rightarrow \cdots \rightarrow a_{km} \rightarrow a_1,
$$
while rotating them by $\frac{2\pi}{m}$.
The mapping class $R_m^{(k)}$ has the property that $(R_m^{(k)})^m$ rotates the interiors of each of the disks by 360 degrees in the counter-clockwise direction.  In other words, it is isotopic
to the product of left Dehn twists along boundary parallel curves.
Thus we have the following.

  \begin{figure}[htbp] 
   \centering
   \includegraphics[width=3in]{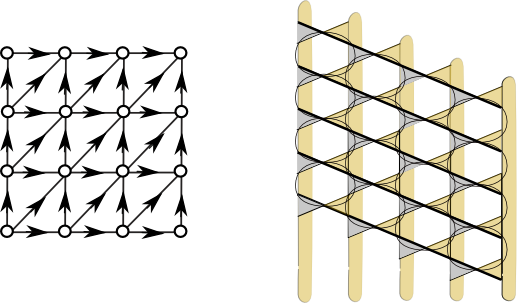}  
   \caption{Seifert surface for full twist briad on $5$ strands and corresponding Coxeter graph $(T_5^{(1)},-1)$.}
   \label{twistgraphsurfaces-fig}
\end{figure}

\begin{lemma}  The twist mapping class $(\Sigma_m^{(k)},R_m^{(k)})$ satisfies
$$
(R_m^{(k)})^{km} = ( \partial_1 \circ \partial_2 \circ \cdots \partial_m)^{-1}
$$
where $\partial_i$ is a positive Dehn twist around the $i$th boundary component of $\Sigma_m$.  (Here, the ordering 
of the $\partial_i$ does
not matter, since the Dehn twists on the right hand side of the equation commute.)
\end{lemma}

We now show that the mapping classes $(\Sigma_m^{(k)},R_m^{(k)})$ are
mixed-sign
Coxeter mapping classes.

Consider the graph $T_m^{(k)}$ with $(km-1)\times (m-1)$ vertices shown in the left diagram of Figure~\ref{twistgraphsurfaces-fig}.  
This graph captures the combinatorics of the set of loops drawn on the surface $\Sigma_m^{(k)}$ on the right of Figure~\ref{twistgraphsurfaces-fig}.  Order the columns of the graph $T_m^{(k)}$ from left to right, and each column from bottom to top.  Starting from the bottom left
corner, moving up the first column, then going to the bottom of the second column to the top, etc.   This ordering is
consistent with arrows in the directed graph shown in  Figure~\ref{twistgraphsurfaces-fig}.  Assign the label $-1$ to
all vertices and denote the signed graph by $(T_m^{(k)},-1)$.   

  \begin{figure}[htbp] 
   \centering
   \includegraphics[width=3.5in]{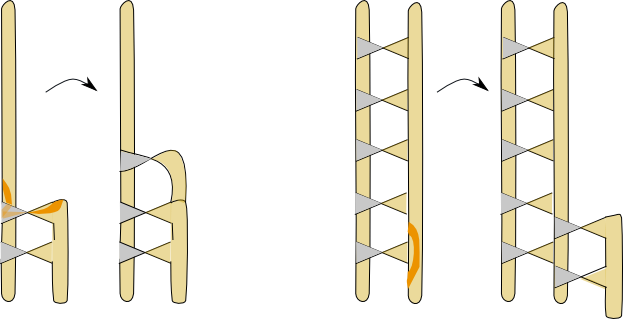}  
   \caption{Constructing the embedding of $\Sigma_m^{(k)}$ in $S^3$ by Hopf plumbings on a disk.}
   \label{TwistHopfPlumb-fig}
\end{figure}

\begin{proposition} The monodromy of $S^3 \setminus L_{m,k}$ with fiber $\Sigma_m^{(k)}$ is the mixed-sign Coxeter mapping class
associated to the labeled graph $(T_m^{(k)},-1)$.
\end{proposition}

\begin{figure}[htbp] 
   \centering
   \includegraphics[width=4.5in]{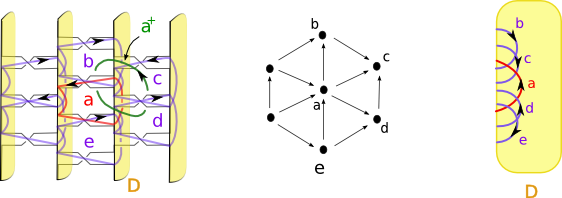}  
   \caption{Local picture of $(\Sigma_5^{(1)},R_m^{(1)})$ .}
   \label{localmesh-fig}
\end{figure}

 \bold{Proof.}   
The surface $\Sigma_m^{(k)} \subset S^3$ is
isotopic to the surface  obtained by successive Hopf plumbing   (see Figure~\ref{TwistHopfPlumb-fig}).
This sequence of clockwise Hopf plumbings is compatible with the ordering on $T_m^{(k)}$.   Furthemore, as seen in 
Figure~\ref{localmesh-fig}, the Seifert matrix is equal to $-I - A^+(T_m^{(k)})$, where $A^+(T_m^{(k)})$ is the upper triangular
part of the adjacency matrix for $T_m^{(k)}$, i.e., the directed adjacency matrix for $T_m^{(k)}$ with the given ordering
of vertices.  The claim thus follows from Proposition~\ref{Seifert-matrix-prop} and Proposition~\ref{Howlett-prop}.
\qed

 \begin{remark}{\em A Coxeter system is {\it spherical} if its associated bilinear form is positive or negative definite, and it is {\it affine}
 if its bilinear form is positive or negative semi-definite.   In the classical case, where $\s \equiv 1$, the Coxeter system is spherical
 or affine if and only if a Coxeter element has spectral radius equal to 1.    In the mixed-sign case, the classification problem 
 is more subtle. 
The graphs $(T_m, -1)$ are examples of 
 mixed-sign Coxeter systems whose Coxeter element has spectral radius one, but is not spherical or affine (for example, 
 the Coxeter group contains hyperbolic elements).   
}
 \end{remark}

\subsection{Flat structure.}

We digress in this section by noting that the mapping classes $(\Sigma_m^{(k)}, R_m^{(k)})$ are naturally endowed with a flat structure,
Furthermore, $(\Sigma_m{(1)},R_m^{(1)})$ are translation surfaces.
Give each $d_1,\dots, d_m$ and $a_1,\dots, a_m$ the flat structure of regular $m$-gons.  Then
$R_m$ preserves the induced singular flat structure on $\Sigma_m$.    More precisely, the $m$-gons $d_1,\dots,d_m$, and $a_1,\dots,a_m$ are 
each permuted cyclically, preserving centers and rotating the $m$-gons by an angle of $\frac{\pi}{m}$.  

  \begin{figure}[htbp] 
   \centering
   \includegraphics[width=3.5in]{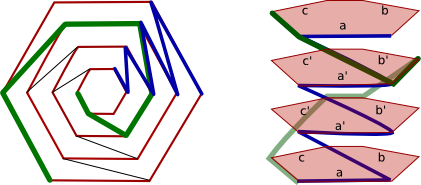}  
   \caption{The main disks $d_1,\dots,d_m, d_{m+1}$, where $d_{m+1}$ is identified with $d_1$ ($m=3$).}
   \label{twistsurface-fig}
\end{figure}

Figure~\ref{twistsurface-fig} shows two views of $\Sigma_m$ for $m=3$.  The left-most figure shows a view of the surface with boundary, where
one of the three boundary curves is drawn spiraling inward.  
The $3$-gons $d_1,d_2,d_3$ are drawn as hexagons (the innermost hexagon is identified with the outermost hexagon), and the boundary
of one of the hexagons $a_i$ is drawn as a zigzag.  The right-hand figure gives a side view, where the top hexagon corresponds
to the inner hexagon in the left hand diagram.   Again the bottom hexagon is
identified with the top hexagon.   The map $R_m$ takes each $d_i$ and $a_i$ and rotates by an angle of $\frac{\pi}{3}$ in the counter-clockwise direction.

  \begin{figure}[htbp] 
   \centering
   \includegraphics[width=3.5in]{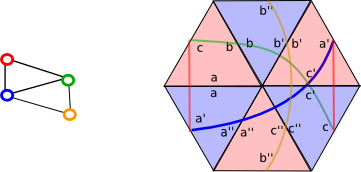}  
   \caption{Two views of the twist surface (homeomorphic to a torus) for $m=3$.}
   \label{flatsurface3-fig}
\end{figure}

If we shrink the boundary of $\Sigma_m$ to a point, then the resulting surface can be given a flat structure as the union of $2m$ regular $m$-gons of
equal size. 
To visualize the flat structure on the
closure $\overline\Sigma_m$, one contracts the spiraling boundary curves.  For example, $\overline \Sigma_3$ is a torus, and $R_3$ preserves
its structure as a union of six equilateral triangles.   In Figure~\ref{flatsurface3-fig} the surface $\overline\Sigma_3$ and its flat structure
are shown.  Sides labeled with the same symbol are identified.  One can see that $\overline\Sigma_3$ is a translation surface.

  \begin{figure}[htbp] 
   \centering
   \includegraphics[width=2.5in]{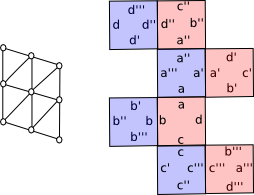}  
   \caption{The twist surface for $m=4$.}
   \label{flatsurface4-fig}
\end{figure}

The map $R_4$ preserves the structure of $\overline \Sigma_4$ as the union of 8 squares (Figure~\ref{flatsurface4-fig}).   Again, we see
that $\overline\Sigma_4$ is a translation surface.  While the flat structure on $\overline \Sigma_3$ has no singularities, the flat
structure on $\overline\Sigma_4$ has $4$ singularities of degree 2.

\subsection{Iterated Murasugi sum with twist maps}   
In this section we use twist maps $(\Sigma_m^{(k)},R_m^{(k)})$ to build 
families of mapping classes with the same mapping torus.
  \begin{figure}[htbp] 
   \centering
   \includegraphics[height=1in]{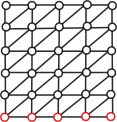}  
   \caption{The extended twist graph  for $m=6$.}
   \label{twistgraph-fig}
\end{figure}

On the level of Coxeter graphs, the construction goes as follows.
Let $(\Gamma,\s)$ be a mixed-sign Coxeter mapping class containing 
$(A_m,-1)$ as a subgraph, where $A_m$ is the standard spherical
Coxeter graph.   We define an extended $m$-twist graph
$(\widetilde T_m,-1)$ with a distinguished $(A_m,-1)$-subgraph, and
define a sequence of graphs $(\Gamma_k,\s_k)$ obtained by iteratively
joining the extended $m$-twist graphs.

\begin{figure}[htbp] 
   \centering
   \includegraphics[height=3in]{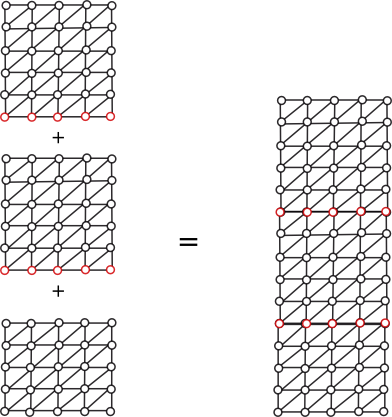}  
   \caption{The iterated twist graph (right) for $m=4$ and $k=3$.}
   \label{iteratedtwists-fig}
\end{figure}

The extended $m$-twist graph is the graph with $m(m-1)$ vertices
shown in Figure~\ref{twistgraph-fig} for $m=6$.   Figure~\ref{iteratedtwists-fig} 
illustrates the $k$ times iterated $m$-twist graph
$(T_m^{k},-1)$ for $k=3$.  This is obtained by identifying the top row of the $m$-twist graph 
with the bottom row of the $m$-twist graph, and repeating.

  Now consider a general mixed-sign Coxeter
graph $(\Gamma,\s)$ containing an $(A_m,-1)$ subgraph.   Let $(S_k,\phi_k)$ be the mapping
class obtained by joining $(\Gamma,\s)$ to the extended $(T_m^{(k)},-1)$ along $(A_m,-1)$ as in Figure~\ref{iteratedtwists-fig}.

\begin{lemma}\label{unif-lem}  The mapping tori  for $(S_k,\phi_k)$, $k \ge 1$, are contained in a single homeomorphism class $M$ of $3$-manifolds,
and the mapping tori $M_k$ for $(S_k,\phi_k)$ are Dehn fillings of $M$ along the tubular neighborhood $N$
of the suspension of $\sO$, with slope $[1:k]$ with respect to some fixed choice of generators of $\pi_1(N)$.
\end{lemma}

\bold{Proof.} This follows from the fact that the mapping tori for $(\Sigma_m^{(k)},R_m^{(k)})$ belong to a single homeomorphism class.
\qed

\begin{remark}{\em The proof above also applies to the iterated  Murasugi sum of $(\Sigma_m^{(k)},R_m^{(k)})$ with an
 arbitrary mapping class $(S,\phi)$ containing an
embedded $P_{2m}$ whose alternating sides are contained in the boundary of $S$. 
}
\end{remark}

\subsection{Attaching tails, and iterated Hopf plumbing.}

Consider the special case when $m=1$.  Then the corresponding Coxeter graph is $(A_k,-1)$, where $A_k$ is
the classical spherical Coxeter system.
In this case, since $A_k$ is a tree, and hence bipartite (see Section~\ref{bipartite-sec}) the signs on the vertices (in this case ``-1") can be replaced by $1$, and we can
take $T_2^{(k)}$ to be the classical $A_k$ diagram.

Let $(S,\phi)$ is any mapping class with an attaching square, and $(S_k,\phi_k)$ is the
sequence of mapping classes obtained by attaching $(\Sigma_1^{(k)},R_1^{(k)})$.   We call 
$(S_k,\phi_k)$ the sequence of mapping classes obtained from $(S,\phi)$ by {\it attaching a tail}.

A {\it Salem-Boyd} sequence of polynomials is a sequence of the form
$$
P_k(x) = x^k Q(x) + Q^*(x)
$$
where $Q(x)$ is a monic integer polynomial, and $Q^*(x) = x^{\deg(Q)} Q(1/x)$ is the 
{\it reciprocal} of $Q(x)$.   These sequences were used in \cite{Salem44} and \cite{Boyd77}
to study properties of Salem numbers.  The house $|P_k|$ of $P_k$ has the following
properties (see,  \cite{Hironaka:Coxeter} Theorem 12).

\begin{theorem}\label{Salem-Boyd-thm} If $P_k(x) = x^k Q(x) + Q^*(x)$ is a Salem-Boyd sequence, then
\begin{description}
\item{(i)} the number of roots of $P_k$ outside the unit circle is monotone increasing, and eventually constant, and
\item{(ii)} $\lim_{k \rightarrow \infty} |P_k| = |Q_k|$.
\end{description}
\end{theorem}

Mapping classes $(S_k,\phi_k)$ corresponding to graphs with tails are studied in \cite{Hironaka:Coxeter}
in the case when the graph is dual to a chord system on a disk.  In this case, it is shown
that the mapping classes are the monodromy of a sequence of links $L_k$ in $S^3$, 
obtained from a single fibered link $L_0$ by twisting a suitable pair of strands.   Furthermore,
the Alexander polynomial, or characteristic polynomial of the action of $\phi_k$ on
first homology, is a Salem-Boyd sequence.  The proof (see \cite{Hironaka:Coxeter}, Theorem 9)
relies only on the form of $(\phi_k)_*$ given in  Proposition~\ref{Howlett-prop} and Proposition~\ref{hom-prop}.

\begin{theorem}[\cite{Hironaka:Coxeter}, Theorem 9]\label{tails-thm}
The Alexander polynomial $\Delta_k$ corresponding to $(S_k,\phi_k)$ is a Salem-Boyd sequence,
and hence the homological dilatations $\lambda_{hom}(\phi_k)$ form a convergent sequence.
\end{theorem}

\begin{corollary}\label{Salem-Boyd-cor1} Let $(S_k,\phi_k)$ be obtained by attaching a single
tail to a mapping class $(S,\phi)$.
Then either 
\begin{description}
\item{(i)} $\lambda_{\mbox{hom}}(\phi_k) = 0$ for all $k$, or
\item{(ii)} $|\Delta_k|$ converges to a real number greater than one.
\end{description}
\end{corollary}

\bold{Proof.} This follows from Theorem~\ref{tails-thm} and Theorem~\ref{Salem-Boyd-thm}.
\qed

\begin{corollary} Let $(\Gamma,\s)$ be a  connected mixed-sign Coxeter graph, $v\in \sV$  a distinguished vertex,
and $(\Gamma_k,\s_k)$ the join of $\Gamma,\s$ with $(A_k,\s(v))$.    Let $(S_k,\phi_k)$ be any
Coxeter mapping class associated to $(\Gamma_k,\s_k)$.    Then either
\begin{description}
\item{(i)} $\phi_k$ is periodic, 
\item{(ii)} $\phi_k$ is pseudo-Anosov, but $\lambda_{\mbox{hom}}(\phi_k) = 0$
for all $k$, or
\item{(iii)} for large $k$, $\phi_k$ is pseudo-Anosov, and there is a constant $C$ such that
$$
\lambda(\phi_k) \ge C > 1.
$$
\end{description}
\end{corollary}

\bold{Proof.} Since $\Gamma$ is connected,  each $\phi_k$ is either periodic, or
$\phi_k$ is pseudo-Anosov.  If $|\Delta_k| > 1$ for some $k$, then, by Corollary~\ref{Salem-Boyd-cor1},
$\lambda_{\mbox{hom}}(\phi_k) =| \Delta_k|$ is greater than one for $k$ large enough.
In this case, by Proposition~\ref{pAcrit-prop},
$\phi_k$ is pseudo-Anosov, and, for any $\epsilon > 0$, $\lambda(\phi_k) \ge \ell - \epsilon$,
where $\ell = \lim_{k \rightarrow\infty}|\Delta_k|$.\qed

 \begin{example}\label{bamboo-ex}{\em
Let $(\Gamma_{m,n},\s)$ be the graph in Figure~\ref{bamboo-fig}.    Since  $T_1^{(k)}$ is a bipartite graph, and hence (by Theorem~\ref{bipartite-thm})  
 it is interchangeable with the standard $A_n$-Coxeter graph with all signs positive.
 Thus, we can think of this graph as being
 obtained from the connected two vertex graph with opposite sign labels on the vertices by joining $m$ and $n$ iterated twist graphs
 of width $r=1$.   
 Since the graph is bipartite, the ordering of the vertices does not change value of $|\omega_{\Gamma_{m,n},\s}|$.
 Let $(S_{m,n},\phi_{m,n})$ be the mapping class associated to $(\Gamma_{m,n},\s)$.  The surface $S_{m,n}$
 has genus
 $$
 g_{m,n} =\left [ \frac{m+n}{2} \right ],
 $$
and one or two boundary components, according to whether $m+n$ is even or odd.
Here  $\left [ a \right ]$ denotes the greatest integer less than or equal to a real number $a$.
In particular, $g_{m,m} = m$.

  \begin{figure}[htbp] 
   \centering
   \includegraphics[width=5in]{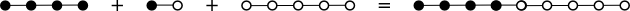}  
   \caption{mixed-sign Coxeter graph $(\Gamma_{4,5},\s)$  obtained by joining twist graphs of width 1.}
   \label{bamboo-fig}
\end{figure}

The mapping classes $(S_{m,n},\phi_{m,n})$ have also been studied in a different form by 
 P. Brinkmann \cite{Brinkmann04} and \cite{HK:braidbounds} (cf. \cite{Tsai08}), yielding the following.
 
\begin{theorem}[Brinkmann \cite{Brinkmann04}, Hironaka-Kin \cite{HK:braidbounds}, Tsai \cite{Tsai08}] For all $m,n$, $(S_{m,n},\phi_{m,n})$ is pseudo-Anosov, and 
for fixed $m+n$, the dilatation is minimized when $m=n$.  Furthermore, 
$$
\log(\lambda(\phi_{g,g})) \asymp \frac{\log(g)}{g}.
$$
\end{theorem}

}
\end{example}

The following question is open.

\begin{question}  Is there a mixed-sign Coxeter graph $(\Gamma,\s)$ with vertices $v_1,\dots,v_m $ (possibly counted with multiplicity) so that
the mapping classes $(S_{\overline k},\phi_{\overline k})$ associated to the joins of $(\Gamma,\s)$ with iterated $k_i$- twists of width 1 at each of the $v_i$,
$\overline k = (k_1,\dots,k_m)$,
 has the asymptotic behavior
\begin{eqnarray}\label{asymp-eqn}
\log(\lambda(\phi_{\overline k})) \asymp \frac{1}{g_{\overline k}},
\end{eqnarray}
where
and $g_{\overline k}$ is the genus of $S_{\overline k}$?
\end{question}

  \begin{figure}[htbp] 
   \centering
   \includegraphics[width=3in]{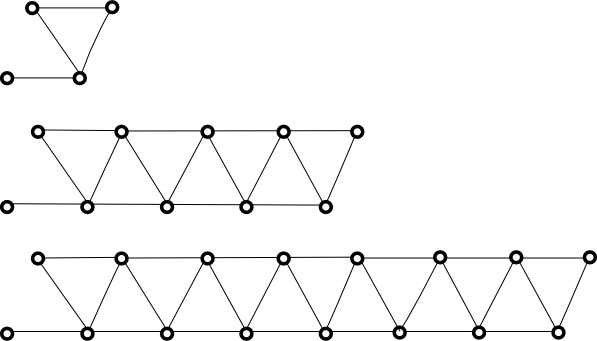}  
   \caption{A negatively signed Coxeter fatgraph for $k=1$, $k=2$ and $k=3$.}
   \label{asympgraph-fig}
\end{figure}

\subsection{Asymptotically small dilatation Coxeter mapping classes.}\label{twistexample-sec}  In this section, 
we prove Theorem~\ref{main2-thm} by showing the existence of a sequence of
 mixed-sign Coxeter mapping classes $(S_k,\phi_k)$ whose closures
 $(\overline S_k,\overline\phi_k)$ have the following properties:
\begin{description}
\item{(i)} $\overline \phi_k$ is pseudo-Anosov, 
\item{(ii)} the associated stable and unstable foliations are orientable,
\item{(iii)} $\overline S_k$ has monotone increasing genus $g_k$,
\item{(iv)} $\lambda(\overline \phi_k)$ approaches 1, and furthermore
$$
\lim_{k\rightarrow \infty} \lambda(\overline \phi_k^{g_k})  = \frac{3 +\sqrt{5}}{2},
$$
the smallest known accumulation point of genus-normalized dilatations.
\end{description}

The mapping classes $(S_k,\phi_k)$ are obtained from
the mixed-sign Coxeter mapping classes corresponding
 to the negatively
signed graphs $\Gamma_k$  drawn in Figure~\ref{asympgraph-fig} by closing over all but 3 boundary
components.   The graphs are
given their fatgraph structures as planar graphs.

\begin{figure}[htbp] 
   \centering
   \includegraphics[height=2in]{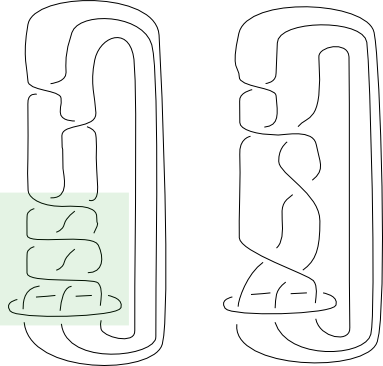}  \hspace{1in}
      \includegraphics[height=2in]{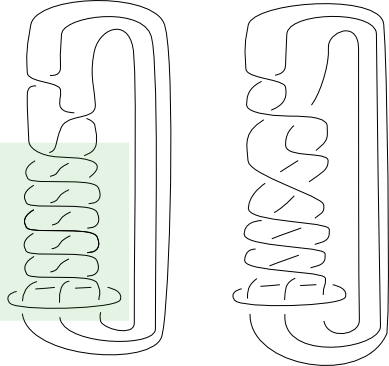}  
   \caption{Link diagram for $S_k,\phi_k$, for $k=1,2$.}
   \label{622multi1-fig}
\end{figure}

To show that  the mapping classes $(S_k,\phi_k)$ satisfy the conditions
(i)-(iv), we relate them to the monodromy of the links drawn
in \ref{622multi1-fig}.  The figure shows pairs of equivalent link diagrams for the
mapping tori of $(S_k,\phi_k)$, $k=1,2$.   In the left versions the shaded region
shows the contribution of the copies of $(\Sigma_3,R_3)$.  The  links beside them on the
right shows an equivalent positive braid version of the links.  The graphs $\Gamma_k$ give
rise to the Seifert surfaces corresponding to the latter planar projection of the links after 
closing over suitable boundary components.


\begin{lemma}\label{example-lem} The mapping classes $(S_k,\phi_k)$ corresponding to $(\Gamma_k,-1)$ satisfy (i)-(iv). \end{lemma}

\begin{figure}[htbp] 
   \centering
      \includegraphics[height=1.5in]{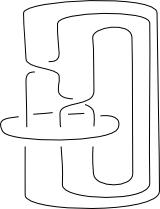}  
   \caption{The $6_{2}^2$-link.}
   \label{622-fig}
\end{figure}

By Lemma~\ref{unif-lem}, the mapping tori for $(S_k^0,\phi_k^0)$ belong to a single homeomorphism class, in this case,
it is the complement $M$ of the $6_2^2$-link $L$ in $S^3$ (see, the knot table in \cite{Rolfsen76}).   The link is shown in 
Figure~\ref{622-fig}.

 The mapping class, written as $\phi_0 = \sigma_1\sigma_2^{-1}$ in terms
of the standard braid generators, is known as the simplest hyperbolic
braid monodromy,   The  fibered face defined by the braid monodromy $\phi_0$
shown in the link diagram in Figure~\ref{622-fig} was studied 
in \cite{Hironaka:LT}.  The link $L$ has two components,  $K_1$ corresponds to the strands of the braid,
and  $K_2$ is the encircling link.   Let $t \in H_1(M,\Z)$ be the meridian loop
around $K_1$ and let $u \in H_1(M;\Z)$ be the meridian loop around $K_2$.   Let $\xi$ be the dual to $t$
and $\psi$ the dual to $u$ in $H^1(M;\Z)$.  Then $\psi$ is the fibration of $M$ corresponding to the
braid monodromy $\phi_0$.   
Given a fibration $\alpha: M \rightarrow S^1$, there is a corresponding element of $H^1(M;\Z)$ defined by
the induced map $\alpha_*: H_1(M;\Z) \rightarrow \Z$, and hence we can write
$$
\alpha_* = a\psi+b\xi.
$$
The integers $a$ and $b$ are determined by the condition that 
$\alpha_*$ restricted to $H_1(S_\alpha;\Z)$ is trivial, where $S_\alpha$ is any fiber of $\alpha$.

For each $k$, we will show that $(S_k,\phi_k)$ is the monodromy of $6k\psi+\xi$.
The link diagrams for $(S_k,\phi_k)$ can be obtained from that of the $6^2_2$-link shown
in Figure~\ref{622-fig}, by removing a tubular neighborhood of $K_2$, and refilling with
meridian $\mu_2'$, where
$$
\mu_2' = \mu_2 -  k\ell_2.
$$
Under the induced homomorphisms
$
\psi_*, \xi_*: H_1(M;\Z) \rightarrow H_1(S^1;\Z) = \Z,
$, we have
$\psi_*(\mu_2') = 1$, and $\xi_*(\mu_2') = -3k$.
Thus, $(S_k,\phi_k)$ is the monodromy of the fibration defined by $\psi_k = 3k \psi + \xi$.

The homological and geometric dilatations of $(S_k,\phi_k)$ can be computed from
the Teichm\"uller and Alexander polynomials of $M$.   
Since $6^2_2$ is a symmetric braid (that is, if one can move the link isotopically to
get the same link diagram with $K_1$ and $K_2$ switched) it follows that the fibered face
for $\psi$ and invariants like the Teichm\"uller polynomial are the same as those of
 the fibered face for $\xi$ (after switching variables).  For the $6^2_2$-link complement,
 the Thurson norm is defined by 
 $$
 ||(a,b)|| = \max\{|2a|,|2b|\},
 $$
 so for $k > 1$, $\psi_k$ lies in the fibered face containing $\xi$, and the rays through
 $\psi_k$ in $H^1(M;\R)$ converge to the ray through $\xi$.  Using the computations in \cite{McMullen:Poly}
 and  \cite{Hironaka:LT},
 we have the following.

\begin{corollary}
The mapping classes $(S_k,\phi_k)$ have genus $3k-1$.
\end{corollary}

\begin{proposition}
The homological and geometric dilatations of $(S_k,\phi_k)$ are given by
$$
\lambda_{\mbox{geo}}(\phi_k) = \lambda(\phi_k) = |LT_{1,3k}| = |x^{6k} - x^{2k+1}-x^{3k}-x^{3k-1} + 1|,
$$
and
$$
\lambda_{\mbox{hom}}(\phi_k) = |x^{6k}-x^{3k+1}+x^{3k} - x^{3k-1} + 1|.
$$
\end{proposition}

\begin{corollary}\label{LT-cor}
When $k$ is even, then $(S_k,\phi_k)$ is orientable, and attains Lanneau and Thiffeault's conjectural minimum
dilatation for orientable pseudo-Anosov mapping classes on closed surfaces of even genus.
\end{corollary}

\begin{corollary} The mapping class $(S_2,\phi_2)$ is a minimum dilatation orientable mapping class for genus $g=5$..
\end{corollary}

The following was proved (in stronger form) in \cite{Hironaka:LT} using fibered face theory.

\begin{lemma}  
$$
\lim_{k \rightarrow \infty}\lambda(\phi_k)^{3k} = \frac{3+\sqrt{5}}{2}.
$$
\end{lemma}

This completes the proof of Theorem~\ref{main2-thm}.

\bibliographystyle{plain}
\bibliography{../../math}
\end{document}